\documentclass[a4,12pt,twoside]{article}
\usepackage{amsfonts}
\usepackage{amssymb}
\usepackage{amsthm}
\usepackage{amsmath}
\usepackage[dvips]{graphicx,color}
\usepackage{latexsym,bm}

\numberwithin{equation}{section}

\renewcommand{\baselinestretch}{1.1}

\def\R2{\mathbb{R}^2}

\begin{document}           

\begin{flushleft}
\renewcommand{\thefootnote}{\fnsymbol{footnote}}
{\bf\LARGE  Exact Hausdorff and packing measures of linear Cantor
sets with overlaps
}  
\footnotetext{{\bf MR(2000) Subject Classification:}
\hspace{0.2cm}28A78, 28A80 }
 \footnotetext{Supported by NSFC 10901081}
 \vskip0.5cm {\large
 Hua Qiu\hspace{0.1cm} }

 {\raggedright \small  Department of Mathematics,
Nanjing University, Nanjing, 210093, China}
\\ { Email: huaqiu@nju.edu.cn}

\end{flushleft}

\noindent{\bf Abstract.}

Let $K$ be the attractor of a linear iterated function system (IFS)
$S_j(x)=\rho_jx+b_j,$ $j=1,\cdots,m$, on the real line satisfying
the generalized finite type condition (whose invariant open set
$\mathcal{O}$ is an interval) with an irreducible weighted incidence
matrix. This condition was introduced by Lau \& Ngai recently as a
natural generalization of the open set condition, allowing us to
include many important overlapping cases. They showed that the
Hausdorff and packing dimensions of $K$ coincide and can be
calculated in terms of the spectral radius of the weighted incidence
matrix. Let $\alpha$ be the dimension of $K$. In this paper, we
state that
\begin{equation*}
\mathcal{H}^\alpha(K\cap J)\leq |J|^\alpha
\end{equation*} for all intervals $J\subset\overline{\mathcal{O}}$, and
\begin{equation*}
\mathcal{P}^\alpha(K\cap J)\geq |J|^\alpha
 \end{equation*} for all intervals $J\subset\overline{\mathcal{O}}$ centered in $K$, where $\mathcal{H}^\alpha$ denotes the $\alpha$-dimensional Hausdorff measure and $\mathcal{P}^\alpha$ denotes the $\alpha$-dimensional packing measure. This result extends a recent work of Olsen where the open set condition is required.
We use these inequalities to obtain some precise density theorems
for the Hausdorff and packing measures of $K$. Moreover, using these
densities theorems, we describe a scheme for computing
$\mathcal{H}^\alpha(K)$ exactly as the minimum of a finite set of
elementary functions of the parameters of the IFS. We also obtain an
exact algorithm for computing $\mathcal{P}^\alpha(K)$ as the maximum
of another finite set of elementary functions of the parameters of
the IFS. These results extend previous ones by Ayer \& Strichartz
and Feng, respectively, and apply to some new classes allowing us to
include linear Cantor sets with overlaps.\hspace{0.2cm}
\\

\section{Introduction and Statement of Results}

In this paper we will analyze the behavior of the Hausdorff and
packing measures of self-similar sets satisfying the generalized
finite type condition which is weaker than the open set condition.
In particular, we will deal with the exact calculating of the
Hausdorff and packing measures for a  kind of linear Cantor sets
with overlaps.

The problem of calculating the  dimension of the attractor of a
self-similar iterated function system (IFS) is one of the most
interesting questions in fractal geometry. During the past two
decades there has been an enormous body of literatures investigating
this problem and wide ranging generalizations thereof. See the books
$\cite{Fal1},\cite{Fal},\cite{Mat1}$ and the references therein. Let
$\{S_j\}_{j=1}^m$ be an IFS of contractive similitudes on
$\mathbb{R}^d$ defined as
\begin{equation}\label{00}
S_j(x)=\rho_jR_j x+b_j,\quad j=1,\cdots,m,
\end{equation}
where  $0<\rho_j<1$ is the contraction ratio, $R_j$ is an orthogonal
transformation and $b_j\in \mathbb{R}^d$, for each $j$. Let $K$
denote the \textit{self-similar set} (or \textit{attractor}) of the
IFS, namely, $K$ is the unique non-empty compact set in
$\mathbb{R}^d$ satisfying
  $$K=\bigcup_{j=1}^mS_j (K).$$
A basic result  (see $\cite{Fal2}$) is that the Hausdorff dimension
$\dim_H K$ and the packing dimension $\dim_P K$ are always equal for
the self-similar set $K$, i.e.,
 $$\dim_H K=\dim_P K.$$
 In general it is quite difficult to calculate this common value, except
 the well-known
classical result (see Moran \cite{Moran}, Hutchinson
\cite{Hutchinson}) that if the IFS satisfies the \textit{open set
condition} (OSC), i.e., there exists a non-empty bounded open set
$\mathcal{O}\subset \mathbb{R}^d$ such that $\bigcup_{j=1}^m
S_j(\mathcal{O})\subset \mathcal{O}$ and $S_i(\mathcal{O})\cap
S_j(\mathcal{O})=\emptyset$ for all $i\neq j$, then the dimension of
$K$ is the unique solution $\alpha$ of the equation
\begin{equation}\label{000}
\sum_{j=1}^m\rho_j^\alpha=1.
\end{equation}
Non-overlapping or almost non-overlapping self-similar IFSs have
been studied in great detail via the OSC.

In the absence of the OSC, much less is known about IFSs with
overlaps.
 To deal with such systems, by extending a method of Lalley $\cite{Lalley}$ and Rao \&Wen $\cite{Rao}$, Ngai \& Wang $\cite{Ngai}$ formulated a weaker separation condition, the \textit{finite type condition} (FTC), which may includes many important overlapping cases  and described a method for computing the Hausdorff and packing dimensions of the attractor in terms of the spectral radius of an associated weighted incidence matrix. The FTC requires the contraction ratios of the IFS's maps to be exponentially commensurable and thus does not generalize the OSC. Recently, Lau and Ngai \cite{Lau} introduced a more general condition, the \textit{generalized finite type condition} (GFTC), which do not need the above requirement, that extends both the OSC and the FTC.
 Under the GFTC, one can also compute the  dimension of the attractor in terms of the spectral radius of a  weighted incidence matrix.

Another central problem concerning  the theory of self-similar IFSs
is to estimate the  Hausdorff and packing measures of self-similar
sets, which is also an area of active research
\cite{Ayer,Feng,Mar,Mar1,Zhou}. In these papers and the references
therein one can find the analysis of the estimation of the values of
the measures of some particular self-similar constructions. Since
the definitions of these measures are sometimes awkward to work
with, there are only very few non-trivial examples of sets whose
exact  measures are known. $\cite{Zhou}$ is a recent review of
relevant open questions in this field. So far most of these
researches have been mainly addressed to the determination of the
upper and lower bounds of the measures. With regard to the
determination of the exact values of them, two papers $\cite{Ayer}$
and $\cite{Feng}$ should be mentioned.

In $\cite{Ayer}$, Ayer and Strichartz considered a kind of linear
Cantor set $K$ which is the attractor of a linear IFS
$S_j(x)=\rho_jx+b_j$, $j=1,\cdots,m$, on the real line satisfying
the OSC (where the open set is the interval (0,1)). Let $\alpha$ be
the  dimension of $K$. They gave an algorithm for computing the
 Hausdorff measure $\mathcal{H}^\alpha(K)$
exactly as the minimum of a finite set of elementary functions of
the parameters of the IFS by using the fact that the exact value of
$\mathcal{H}^\alpha(K)$ is the inverse of the maximal density of
intervals  contained in $[0,1]$ with respect to the
 normalized measure $\lambda$ of
$\mathcal{H}^\alpha$ restricted to $K$, where
$\lambda=\mathcal{H}^\alpha|_K/\mathcal{H}^\alpha(K)$. It should be
pointed out that if the OSC is satisfied, then
$\mathcal{H}^\alpha|_K$ and $\mathcal{P}^\alpha|_K$  are
proportional. Hence  $\lambda$ is also equal to
$\mathcal{P}^\alpha|_K/\mathcal{P}^\alpha(K)$.

On the other hand, in $\cite{Feng}$, Feng proved that the
 packing measure $\mathcal{P}^\alpha(K)$  is
equal to the inverse of the so-called minimal centered density of
intervals centered in $K$ with respect to  $\lambda$, which also
yields an explicit formula for calculating the exact value of
$\mathcal{P}^\alpha(K)$  in terms of the parameters of the IFS.

However, in these papers, one needs to work on self-similar sets
under the OSC. To the best of our knowledge, there is no result
concerning the exact Hausdorff or packing measures of self-similar
sets without the OSC. Since the calculation of the dimension of
self-similar sets under the OSC can be successfully extended to
those sets satisfying the GFTC which includes many interesting
overlapping cases, and in view of the above discussion, it is
natural to ask whether the Hausdorff or packing measure of the
linear Cantor set $K$  under only the GFTC can also be calculated
exactly. This is the main goal of this paper.

In $\cite{Ayer}$ and $\cite{Feng}$,  to get the exact values of
$\mathcal{H}^\alpha(K)$ and $\mathcal{P}^\alpha(K)$ of  $K$ under
the OSC (where the open set is the  interval (0,1)), the following
explicit formulae play a key role.
\begin{equation}\label{aa}
\mathcal{H}^\alpha(K)^{-1}=\sup\{\frac{\lambda(J)}{|J|^\alpha}: J
\mbox{ is an interval with } J\subset [0,1]\},
\end{equation}
and
\begin{equation}\label{bb}
\mathcal{P}^\alpha(K)^{-1}=\inf\{\frac{\lambda(J)}{|J|^\alpha}: J
\mbox{ is an interval centered in } K \mbox{ with } J\subset
[0,1]\}.
\end{equation}
 Formula $(\ref{aa})$ was
implicit in earlier work by Marion $\cite{Mar,Mar1}$ and in
$\cite{Ayer}$ by Ayer \& Strichartz,  while formula $(\ref{bb})$ was
proved in a direct and elementary way in $\cite{Feng}$ by Feng.

Recently, Mor\'{a}n $\cite{Mo}$ and Olsen $\cite{Ols1}$ extended the
above two formulae to the higher dimensional case independently. In
$\cite{Mo}$ the so-called self-similar tiling principle plays a
central role in the proof. This principle says that any open subset
$U$ of $K$ can be tiled by a countable set of similar copies of an
arbitrarily given closed set with positive Hausdorff or packing
measure while the tiling is exact in the sense that the part of $U$
which cannot be covered by the tiles is of null measure.

The proof in $\cite{Ols1}$ is quite different from that in
$\cite{Mo}$. Let $K$ be a self-similar set in $\mathbb{R}^d$ as
described in $(\ref{00})$ with the  dimension $\alpha$, under the
OSC or the \emph{strong separation condition} (SSC). Recall that in
$\cite{Ols1}$ Olsen performs a detailed analysis of the behavior of
the Hausdorff measure $\mathcal{H}^\alpha(K\cap U)$ and  the packing
measure $\mathcal{P}^\alpha(K\cap B(x,r))$ of small convex Borel
sets $U$ and balls $B(x,r)$. In particular, he showed that if $K$ is
under the OSC, then
$$\mathcal{H}^\alpha(K\cap U)\leq |U|^\alpha$$ for each convex Borel set $U$. A dual result for the packing measure was also proved which says that
$$\mathcal{P}^\alpha(K\cap B(x,r))\geq(2r)^\alpha$$ for each $x\in K$ and small $r>0$ if $K$ satisfies the SSC. The latter result was generalized to the OSC case
 by the author recently (see \cite{Qiu}) in proving the continuity of the packing measure function of self-similar IFSs, which says that  the above inequality  actually holds for each $B(x,r)$ contained in $\mathcal{O}$ with $x\in K$, where $\mathcal{O}$
  is an open set associated with the OSC, satisfying $\mathcal{O}\cap K\neq\emptyset$. (There must exist such $\mathcal{O}$ since the OSC is equivalent to the strong OSC. See $\cite{Shi}$.)

To  match our question precisely we restrict our interest to the
linear Cantor set $K$ defined before. Hence the above two formulae
are rewroted in the following form, i.e.,
\begin{equation}\label{cc}
\mathcal{H}^\alpha(K\cap J)\leq |J|^\alpha
\end{equation} for all intervals $J\subset [0,1]$, and
\begin{equation}\label{dd}
\mathcal{P}^\alpha(K\cap J)\geq |J|^\alpha
 \end{equation}
 for all intervals $J\subset [0,1]$ centered in $K$.

As an application of $(\ref{cc})$ and $(\ref{dd})$, Olsen reproved
formulae $(\ref{aa})$ and $(\ref{bb})$ (he proved the general higher
dimensional case) using the classical density theorems of geometric
measure theory which were stated for arbitrary subsets of Euclidean
space $\cite{Fal1,Mat1}$.

 Formulae $(\ref{aa})$ and
$(\ref{bb})$ say that the exact values of $\mathcal{H}^\alpha(K)$
and $\mathcal{P}^\alpha(K)$ coincide with the inverses of the
supremum and infimum of the densities of $\lambda$ on suitable
classes of sets, respectively.

In order to calculate the exact values of the measures of  $K$
satisfying  the GFTC, we need to establish analogous explicit
formulae of $(\ref{aa})$ and $(\ref{bb})$. Following the frame of
Olsen's work, two inequalities similar to $(\ref{cc})$ and
$(\ref{dd})$ are required. Recall that in proving $(\ref{cc})$ and
$(\ref{dd})$, one should find optimal coverings and packings in a
self-similar setting which require almost non-overlap among the
various similar pieces into which the fractal decomposes. In view of
this, it is therefore entirely plausible that the OSC is
indispensible. In the present paper, somewhat surprisingly, we will
show that the formulae $(\ref{cc})$ and $(\ref{dd})$  still hold
under the assumption that the weighted incidence matrix of $K$ is
irreducible where $K$ is required to satisfy only the GFTC. This
leads to the following results.

\textbf{Theorem 1.1. } \emph{Let $K$ be a linear Cantor set
satisfying the GFTC with respect to the invariant open set $(0,1)$
with an irreducible weighted incidence matrix $A_\alpha$, where
$\alpha$ is the Hausdorff dimension of $K$. Then}
\begin{equation}\label{2}
\mathcal{H}^\alpha(K\cap J)\leq |J|^\alpha
\end{equation} \emph{for all intervals $J\subset [0,1]$.}

\textbf{Theorem 1.2. } \emph{Let $K$ be the linear Cantor set described as before. Then
\begin{equation}\label{9}
\mathcal{P}^\alpha(K\cap J)\geq |J|^\alpha
\end{equation} for all intervals $J\subset [0,1]$ centered in $K$.}

We will give the detailed definition of  the GFTC and the exact
concept of the weighted incidence matrix $A_\alpha$ in Section 2.

The idea of establishing Theorem 1.1 and Theorem 1.2 is the
following. We first observe that if the weighted incidence matrix
$A_\alpha$ of $K$ is irreducible, then $K$ can be decomposed into an
union of a set $K_a$ with a graph directed construction  and an
attractor $K_b$ of a countable infinite IFS under the OSC (see
\cite{Mau}, \cite{Mau1} for further properties of the graph directed
sets and the infinite IFSs, respectively). Moreover, the  dimension
of $K_a$ is strictly less than that of $K_b$. Hence the subset $K_a$
will have null $\alpha$-dimensional Hausdorff (packing) measure
which ensures us to consider $K_b$ in stead of $K$. Noticing that
$K_b$ is an attractor of a countable infinite IFS satisfying the
OSC, it is possible to adapt the techniques for proving $(\ref{cc})$
and $(\ref{dd})$ to establish Theorem 1.1 and Theorem 1.2.

 By a similar discussion  for self-similar sets under the OSC, the $\alpha$-dimensional Hausdorff  measure restricted to $K_b$ and the $\alpha$-dimensional packing measure restricted
  to $K_b$ are also proportional. Obviously the above fact still holds if we replace $K_b$ by $K$. We still write the  normalized  measure of $\mathcal{H}^\alpha$ restricted to $K$ as $\lambda$, then
$\lambda=\mathcal{H}^\alpha|_K/\mathcal{H}^\alpha(K)=\mathcal{P}^\alpha|_K/\mathcal{P}^\alpha(K)$.
We will show the $\lambda$-measure of some  special
 kind of sets called \emph{islands} of $K$ can be expressed in terms of the parameters of the IFS of $K$.
  Then following the frame of Olsen's work, we use the inequalities $(\ref{2})$ and $(\ref{9})$ to get the following explicit formulae for $\mathcal{H}^\alpha(K)$ and $\mathcal{P}^\alpha(K)$, analogous to $(\ref{aa})$ and $(\ref{bb})$.

\textbf{Corollary 1.3.}  \emph{Let $K$ be the linear Cantor set described as before. Then
\begin{equation}\label{8}
\mathcal{H}^\alpha(K)^{-1}=\sup\{\frac{\lambda(J)}{|J|^\alpha}: J \mbox{ is an interval with } J\subset [0,1]\}.
\end{equation}}

\textbf{Corollary 1.4.}
\emph{Let $K$ be the linear Cantor set described as before. Then
\begin{equation}\label{16}
\mathcal{P}^\alpha(K)^{-1}=\inf\{\frac{\lambda(J)}{|J|^\alpha}: J \mbox{ is an interval centered in } K \mbox{ with } J\subset [0,1]\}.
\end{equation}
}

 These corollaries extend the results in $\cite{Mo,Ols1,Qiu}$.
Following the technique frame of $\cite{Ayer,Feng}$, under suitable
assumptions, we then give an algorithm for computing
$\mathcal{H}^\alpha(K)$ and $\mathcal{P}^\alpha(K)$ exactly as the
inverse of the maximal or minimal value of suitable finite sets of
elementary functions of the parameters of the IFS respectively. This
is possible since we could make a detailed analysis of $\lambda$,
and thus a detailed analysis of the supremum in $(\ref{8})$ and the
infimum in $(\ref{16})$ respectively. It should be  mentioned here
that we may allow touching islands, and indeed this case will lead
to some of complicated and interesting phenomena. Due to the fact
that the self-similar construction of $K$ under the GFTC is much
more complicated than that under the OSC, our description of the
exact calculations of the two kinds of measures will need some new
important notations and techniques. We will describe a big scheme
for the exact computing, which is a major adaption of the techniques
used in $\cite{Ayer}$ and $\cite{Feng}$.

This paper is organized as follows. In Section 2, we give some
notations and basic facts about the GFTC. Our description  of the
GFTC is slightly different but equivalent to the original version in
$\cite{Lau}$. In Section 3, we deal with the density theorems for
the Hausdorff and packing measures of linear Cantor sets satisfying
the GFTC. Firstly, we give the proofs of Theorem 1.1 and Theorem 1.2
respectively. Secondly, we prove the formulae in Corollary 1.3 and
Corollary 1.4 using the classical density theorems of geometry
measure theory. Throughout this section and the following ones, we
will always assume that the weighted incidence matrix of $K$ is
irreducible. In Section 4, we focus on the calculation of the exact
measures of the linear Cantor $K$ under some suitable assumptions. A
 scheme is provided.  Section 5 collects some further discussions
on this subject. We consider the possibility of dropping some
assumption required in Section 4. We discuss briefly the slightly
more general cases of IFSs that contain orientation reversing
simlarities. We also consider the situation in higher dimensional
Euclidean spaces and show why our results can not be generalized.
Throughout the context, we will show some interesting and
non-trivial examples.

\section{Linear Cantor sets under the GFTC}

For convenience, we introduce a slightly different but equivalent
description of  the GFTC defined in \cite{Lau}. We will focus the
interest on the linear Cantor sets on the real line $\mathbb{R}$.

We will use the following notations throughout the paper. For a
Borel measure $\nu$ on $\mathbb{R}$ and a Borel set $E$, we let
$\nu|_E$ denote the restriction of $\nu$ to $E$. For any subset
$E\subset\mathbb{R}$, we denote the \emph{diameter} of $E$ by $|E|$.
For any $x\in \mathbb{R}$, let $\mbox{dist}(x,E)$ denote the
distance between $x$ and $E$, namely, $\mbox{dist}(x,E)=\inf\{|x-y|:
y\in E\}.$ If $A$ is any finite or countable set, we denote by
$\sharp A$ the \emph{cardinality} of $A$.  For $E\subset
\mathbb{R}^d$, $s\geq 0$ and $\delta>0$, put
$\mathcal{H}^s_\delta(E):=\inf\{\sum_{i}|U_i|^s\},$ where the
infimum is taken over all \emph{$\delta$-coverings} of $E$, i.e.,
countable collections $\{U_i\}$ of subsets of $\mathbb{R}^d$ with
diameters smaller than $\delta$ such that $E\subset\bigcup_{i}U_i$.
The \emph{$s$-dimensional Hausdorff measure} $\mathcal{H}^s(E)$ of
$E$ is defined by
$$\mathcal{H}^s(E):=\sup_{\delta>0}\mathcal{H}^s_\delta(E).$$
The \textit{Hausdorff dimension} of $E$ is defined as
$$\dim_{H}E:=\inf\{s\geq 0| \mathcal{H}^s(E)=0\}=\sup\{s\geq 0| \mathcal{H}^s(E)=\infty\}.$$
Recall the definition of the packing measure, introduced by Tricot
\cite{Tri}, Taylor \& Tricot \cite{TT}, which requires two limiting
procedures. For $E\subset\mathbb{R}^d$ and $\delta>0$, a
\emph{$\delta$-packing} of $E$ is a countable family of disjoint
open balls of radii at most $\delta$ and with centers in $E$. For
$s\geq 0$, the \textit{$s$-dimensional packing premeasure} of $E$ is
defined as
$$P^s(E):=\inf_{\delta>0}\{P^s_\delta(E)\},$$
where $P^s_\delta(E):=\sup\{\sum_{i}|B_i|^s\}$ with the supremum taken over all $\delta$-packing of $E$.  The \textit{$s$-dimensional packing measure} of $E$ is defined as
$$\mathcal{P}^s(E):=\inf\{\sum_{i}P^s(E_i)| E\subset\bigcup_{i} E_i\}.$$
The \textit{packing dimension} of $E$ is defined as
$$\dim_{P}E:=\inf\{s\geq 0| \mathcal{P}^s(E)=0\}=\sup\{s\geq 0| \mathcal{P}^s(E)=\infty\}.$$  See \cite{Fal} and \cite{Mat1} for further properties of the above measures and dimensions.

 Let $\{S_1,\cdots,S_m\}$ be a linear IFS of contractive similitudes on the line $\mathbb{R}$ defined by
\begin{equation}\label{1}
S_j(x)=\rho_jx+b_j, j=1,\cdots,m,
\end{equation}
with contraction ratios satisfying $0<|\rho_j|<1$. Let $K$ denote
its attractor. Let $\Sigma=\{1,\cdots,m\}$ and let
$\Sigma_*:=\bigcup_{k=0}^\infty\Sigma_k$ be the symbolic space
representing the IFS (by convention, $\Sigma_0=\emptyset$).  For
$\textbf{i}=(i_1,\cdots,i_k)\in \Sigma_*$, $E\subset \mathbb{R}$, we
use the standard notation $S_\textbf{i}:=S_{i_1}\circ\cdots\circ
S_{i_k},  \rho_\textbf{i}:=\rho_{i_1}\cdots\rho_{i_k},
E_{\textbf{i}}=S_\textbf{i}(E),$ with $S_\emptyset:=I,$ the identity
map on $\mathbb{R}$, and $\rho_{\emptyset}:=1$. Let
$\textbf{i}=(i_1,\cdots,i_k)$ and $\textbf{j}=(j_1,\cdots,j_{k'})$
be two indices of $\Sigma_*$. The \textit{length} of $\textbf{i}$ is
$|\textbf{i}|=k$. We write $\textbf{i}\preceq\textbf{\textbf{j}}$ if
$\textbf{i}$ is an initial segment of $\textbf{j}$, and write
$\textbf{i}\npreceq \textbf{j}$ if $\textbf{i}$ is not an initial
segment of $\textbf{j}$. $\textbf{i}$ and $\textbf{j}$ are
\textit{incomparable} if neither $\textbf{i}\preceq \textbf{j}$ nor
$\textbf{j}\preceq \textbf{i}$.

Let $\{\mathcal{M}_k\}_{k=0}^{\infty}$ be a sequence of index sets, where $\mathcal{M}_k\subset\Sigma_*$ for all $k\geq 0$ and $\mathcal{M}_0=\Sigma_0$. We say that $\{\mathcal{M}_k\}_{k=0}^{\infty}$ is a \textit{sequence of nested index sets} if it satisfies the following conditions:

(1) Both $\{\min\{|\textbf{i}|: \textbf{i}\in\mathcal{M}_k\}\}_{k=0}^\infty$ and $\{\max\{|\textbf{i}|: \textbf{i}\in\mathcal{M}_k\}\}_{k=0}^\infty$ are non-decreasing and have infinity limit;

(2) For each $k\geq 0$, all $\textbf{i},\textbf{j}\in \mathcal{M}_k$ are incomparable if $\textbf{i}\neq \textbf{j}$;

(3) For each $\textbf{j}\in \Sigma_*$ with $|\textbf{j}|>\max\{|\textbf{i}|: \textbf{i}\in\mathcal{M}_k\}$, there exists $\textbf{i}\in \mathcal{M}_k$ such that $\textbf{i}\preceq\textbf{j}$;

(4) For each $\textbf{j}\in \Sigma_*$ with $|\textbf{j}|<\min\{|\textbf{i}|: \textbf{i}\in\mathcal{M}_k\}$, there exists $\textbf{i}\in \mathcal{M}_k$ such that $\textbf{j}\preceq\textbf{i}$;

(5) There exists a positive integer $L$ such that for all
$\textbf{i}\in \mathcal{M}_k$ and $\textbf{j}\in \mathcal{M}_{k+1}$
with $\textbf{i}\preceq \textbf{j}$, we have
$|\textbf{j}|-|\textbf{i}|\leq L$, where $L$ is independent of $k$.

For general sequences, we allow
$\mathcal{M}_k\cap\mathcal{M}_{k+1}\neq \emptyset$ and
$\bigcup_{k=0}^{\infty}\mathcal{M}_k$ may be a proper subset of
$\Sigma_*.$ If we let $\mathcal{M}_k=\Sigma_k$ for all $k\geq 0$, we
get a canonical such sequence. For $k\geq 0$, let
$\Lambda_k:=\{\textbf{i}=(i_1,\cdots,i_n)\in\Sigma_*:
|\rho_\textbf{i}|\leq \rho_{\min}^k<|\rho_{i_1,\cdots,i_{n-1}}|\},$
where $\rho_{\min}=\min\{|\rho_j|:1\leq j\leq m\}$ It is easy to
check that $\{\Lambda_k\}_{k=0}^\infty$ is also a sequence of nested
index sets, which is used to define the FTC in \cite{Ngai}.

Fix a sequence of nested index sets
$\{\mathcal{M}_k\}_{k=0}^\infty$. Note that if $\mathcal{O}\subset
\mathbb{R}$ is a non-empty bounded open set which is
\textit{invariant} under $\{S_j\}_{j=1}^m$, i.e., $\bigcup_{j=1}^m
\mathcal{O}_j\subset \mathcal{O}$, then $\{\bigcup_{\textbf{i}\in
\mathcal{M}_k}\mathcal{O}_\textbf{i}\}_{k=0}^\infty$ becomes a
sequence of nested subsets of $\mathbb{R}$. For each integer $k\geq
0$, let $\mathcal{V}_k$ be the set of \textit{vertices} defined as
$\mathcal{V}_k:=\{(S_\textbf{i},k): \textbf{i}\in \mathcal{M}_k\}.$
We call $(I,0)$ the \textit{root vertex} and let
$\mathcal{V}:=\bigcup_{k\geq 0}\mathcal{V}_k$. Note that if
$S_\textbf{i}=S_\textbf{j}$ for some $\textbf{i}\neq \textbf{j}\in
\mathcal{M}_k$, they determine the same vertex. For
$\textbf{v}=(S_\textbf{i},k)\in \mathcal{V}_k$, we introduce the
convenient notation $S_\textbf{v}:=S_\textbf{i}$ and
$\rho_\textbf{v}:=\rho_\textbf{i}$. The notation $S_\textbf{v}$
allows us to refer to a vertex in $\mathcal{V}_k$ without explicitly
specifying the index $\textbf{i}$.

For simplicity, we assume that there is a following form of
invariant set which will be used in the definition of the GFTC: an
open interval $\mathcal{O}$ which is invariant under
$\{S_j\}_{j=1}^m$. (There are examples where the GFTC holds, but not
with an open interval.) Without loss of generality we take
$\mathcal{O}=(0,1)$.

For any $k\geq 0$, let $F_k=\bigcup_{\textbf{v}\in
\mathcal{V}_k}\overline{\mathcal{O}}_\textbf{v}$. Notice that
$\overline{\mathcal{O}}_\textbf{v}\subset \overline{\mathcal{O}}$.
Then we have
$$K=\bigcap_{k=0}^\infty F_k.$$
For each $\textbf{v}\in \mathcal{V}_k$, $\overline{\mathcal{O}}_\textbf{v}$ is a sub-interval contained in $[0,1]$ with endpoints $S_\textbf{v}(0)$ and $S_\textbf{v}(1)$. We call $\overline{\mathcal{O}}_\textbf{v}$ a $k$\textit{-th generation interval} of $K$.

We call two intervals are \textit{separate} if they have at most one common point. Otherwise, we call they are \emph{overlapping}. From the definition of $F_k$, we see that $\bigcup_{\textbf{v}\in \mathcal{V}_k}\mathcal{O}_\textbf{v}$ consists of some separate open intervals, and each open interval is an union of the interiors of one or several $k$-th generation intervals. We call the closure of each such open interval a \emph{$k$-th generation island}, and use $\mathcal{F}^0_k$ to denote the set of all $k$-th generation islands, and $\mathcal{F}_k$ the finite field generated from $\mathcal{F}^0_k$. We call the unique element $\overline{\mathcal{O}}=[0,1]$ in $\mathcal{F}_0^0$ the \emph{root island}. For the open intervals between each pair of the $k$-th generation islands, we call \emph{lakes}. For each $k$-th generation island $I$, we use $V(I)$ to denote the vertices set of all $k$-th generation intervals contained in $I$, i.e.,
$$V(I)=\{\textbf{v}\in \mathcal{V}_k: \overline{\mathcal{O}}_\textbf{v}\subset I\}.$$ It is easy to verify that $I=\bigcup_{\textbf{v}\in V(I)}\overline{\mathcal{O}}_\textbf{v}$. We call each such interval $\overline{\mathcal{O}}_\textbf{v}$ a \emph{constitutive interval} of $I$.
Let $I\in \mathcal{F}^0_k$ and $I'\in \mathcal{F}^0_{k+1}$ for some
$k\geq 0$. Then either $I'\subset I$ or they are separate. If it is
 the first case, we call $I$ a \emph{parent} of $I'$ and $I'$ an
\emph{offspring}(or \emph{descendant}) of $I$.

We define an equivalence relation on  $\mathcal{F}^0:=\bigcup_{k\geq 0}\mathcal{F}^0_k$ to identify islands that are isomorphic in the sense that they behave the same overlap type.

\textbf{Definition 2.1.} \emph{Two islands $I\in \mathcal{F}_k^0$
and $I'\in \mathcal{F}^0_{k'}$ are \emph{equivalent}, denoted by
$I\sim I'$, if there is a linear function $\tau$ mapping $I$ onto
$I'$, such that the following conditions are satisfied:}

\emph{(1) $\{S_{\textbf{v}'}: \textbf{v}'\in V(I')\}=\{\tau\circ
S_\textbf{v}: \textbf{v}\in V(I)\}$;}

\emph{(2) For any positive integer $l\geq 1$, there is an island
$J\in \mathcal{F}_{k+l}^0$ contained in $I$ if and only if there is
also an island $J'\in \mathcal{F}_{k'+l}^0$ contained in $I'$ where
$J'=\tau(J)$.}

It is easy to see that $\sim$ is an equivalence relation. We denote the equivalence class containing $I$ by $[I]$ and call it the \emph{overlap type} of $I$. Condition (2) says that any two islands with the same overlap type have equivalent offsprings.

\textbf{Definition 2.2.} \emph{We say that a linear IFS of
contractive similitudes on $\mathbb{R}$
 satisfies the generalized finite type condition (GTFC), with respect to  the invariant set $\mathcal{O}=(0,1)$, if $\mathcal{O}$ is an invariant set under the IFS and there is a sequence of nested index sets $\{\mathcal{M}_k\}_{k=0}^\infty$, such that $\mathcal{F}^0/\sim=\{[I]: I\in \mathcal{F}^0\}$ is a finite set.
}

It is easy to see that the IFS satisfies the GTFC if and only if
there exists some $k_0\geq 0$ such that none of the islands in
$\mathcal{F}^0_{k_0+1}$ is of a new overlap type. Let
$\mathcal{T}_1,\cdots,\mathcal{T}_q$ denote all the distinct overlap
types, with $\mathcal{T}_1=[\overline{\mathcal{O}}]$.

 For each
$\alpha\geq 0$ we define a \emph{weighted incidence matrix}
$A_\alpha=A_\alpha(i,j)_{i,j=1}^q$ as follows. Fix $i(1\leq i\leq
q)$ and an island $I\in \mathcal{F}^0$ with $[I]=\mathcal{T}_i$.
Suppose that $I$ is a $k$-th generation island, let $I_1,\cdots,I_l$
be the offsprings of $I$ in $\mathcal{F}^0_{k+1}$. Then we define
$$A_\alpha(i,j):=\sum\{\left(\frac{|I_s|}{|I|}\right)^\alpha: [I_s]=\mathcal{T}_j, 1\leq s\leq l\}.$$
It is easy to see that the definition of $A_\alpha(i,j)$ is
independent of the choice of $I$.

\textbf{Theorem 2.3$^{\cite{Lau}}$}. \emph{Let $\lambda_\alpha$ be
the spectral radius of the associated weighted incidence matrix
$A_\alpha$. Then
$$\dim_{H}K=\dim_P K=\alpha,$$
where $\alpha$ is the unique number such that $\lambda_\alpha=1$.
Moreover, $0<\mathcal{H}^\alpha(K)<\infty$.}

For the convenience of the readers, we would like to give a direct
and elementary proof of this result. The proof makes use of the
ideas in $\cite{Lau}$ and $\cite{Mau}$.

\textit{Proof.} We need to define a natural probability measure
$\mu$ on $K$. Let $(a_1,\cdots,a_q)^T$ be an $1$-eigenvector of
$A_\alpha$, normalized so that $a_1=1$. (This is possible because
all overlap types are descendants of $\mathcal{T}_1$.) Here $\alpha$
is the unique number such that $\lambda_\alpha=1$. For each island
$I$, where $I\in \mathcal{F}^0_k$ and $[I]=\mathcal{T}_i$ for some
$k\geq 0$ and $1\leq i\leq q$, we let $\mu(I)=|I|^\alpha a_i.$
Obviously, $\mu([0,1])=a_1=1.$

To show that $\mu$ is indeed a probability measure on $K$, we notice
that two islands $I\in \mathcal{F}^0_k$ and $I'\in
\mathcal{F}^0_{k'}$ with $k\leq k'$, are overlapping if and only if
either $I=I'$ in the case $k=k'$ or $I'$ is a descendant of $I$ in
the case $k<k'$. In both cases, $I'\subset I$. Now let $I\in
\mathcal{F}^0_k$ and let $\mathcal{D}$ denote the set of all
offsprings of $I$ in $\mathcal{F}_{k+1}^0$. Then
\begin{eqnarray*}
\sum_{I'\in \mathcal{D}}\mu(I')&=&\sum_{j=1}^q\sum\{\mu(I'): I'\in
\mathcal{D}, [I']=\mathcal{T}_j\}
= |I|^\alpha \sum_{j=1}^q\sum\{\left(\frac{|I'|}{|I|}\right)^\alpha a_j: I'\in \mathcal{D}, [I']=\mathcal{T}_j\}\\
&=& |I|^\alpha\sum_{j=1}^q A_\alpha(i,j)a_j= |I|^\alpha a_i=\mu(I).
\end{eqnarray*}
It follows now from $\mu([0,1])=1$ that $\mu$ is indeed a
probability measure on $K$.

\emph{Lower bound.} Let $E$ be a bounded Borel subset of
$\mathbb{R}$ and let $\mathcal{N}(E)$ be defined  as
$$\mathcal{N}(E):=\{I\in \mathcal{F}^0: |I|\leq |E|<|\mathcal{P}(I)| \mbox{ and } I\cap E\neq \emptyset\},$$ where $\mathcal{P}(I)$ denotes the parent of the island $I$.
 It is easy to verify that for any bounded Borel set $E\subset \mathbb{R}$, $\sharp \mathcal{N}(E)\leq C_0:=\max\{{|\mathcal{P}(I)|}/{|I|}: I\in \mathcal{F}^0\}+2$. Note that $\mu(E)\leq \sum _{I_j\in \mathcal{N}(E)}\mu(I_j).$
If we assume that $[I_j]=\mathcal{T}_{i_j}$, then
$$\mu(E)\leq\sum_{I_j\in \mathcal{N}(E)}|I_j|^\alpha a_{i_j}\leq |E|^\alpha\sum _{I_j\in \mathcal{N}(E)}a_{i_j}\leq C_0\max_{1\leq i\leq q}a_{i}|E|^\alpha.$$
Thus $\mathcal{H}^\alpha (K)>0$ and $\dim_{H}K\geq \alpha$(see
\cite{Fal}), which is the required lower bound.

\emph{Upper bound.} To obtain the upper bound $\dim_H K\leq \alpha,$
we first assume that $A_\alpha$ is irreducible and thus all the
$a_i$'s are positive. For each $k\geq 0$, $K\subset \bigcup_{I\in
\mathcal{F}^0_k}I$ and
$$
\sum_{I\in
\mathcal{F}^0_k}|I|^\alpha=\sum_{i=1}^q\sum\{\frac{1}{a_i}|I|^\alpha
a_i: I\in \mathcal{F}_{k}^0, [I]=\mathcal{T}_i\}\leq \max_{1\leq
i\leq q}\{\frac{1}{a_i}\}\sum_{I\in
\mathcal{F}^0_k}\mu(I)=\max_{1\leq i\leq q}\{\frac{1}{a_i}\}<\infty.
$$
Since for each $k\geq 0$, $\mathcal{F}^0_k$ is a covering of $K$,
and $\lim_{k\rightarrow\infty}\max\{|I|: I\in \mathcal{F}^0_k\}=0$,
the definition of the Hausdorff measure implies that
$\mathcal{H}^\alpha(K)<\infty$, and thus $\dim_H K\leq \alpha$.

Now assume $A_\alpha$ is not irreducible. After an appropriate
permutation we can assume that $A_\alpha$ has the form
$$A_\alpha=\left[\begin{array}{cccc}
               A_1 & * & \cdots & * \\
                0 & * & \cdots & * \\
                \vdots &  & \cdots & * \\
                0 &\cdots & \cdots & A_r \\
                        \end{array}\right],$$
where each $A_i$ is either an irreducible square matrix or an
$1\times 1$ zero matrix. Let $\mathcal{E}:=\{A_i: A_i  \mbox{ is
non-zero}\}$, counting multiplicity. For each $A_i$, let
$\mathcal{T}_{A_i}$ be the collection of overlap types corresponding
to $A_i$. For each $A_i\in\mathcal{ E}$ and each island $I$ with
$[I]\in \mathcal{T}_{A_i}$, define a subset $K_{A_i}(I)\subset K$ as
follows.
$$K_{A_i}(I):=\bigcap_{k'=0}^\infty \bigcup\{I'\in \mathcal{F}^0_{k+k'}: [I'],[\mathcal{P}(I')],\cdots, [\mathcal{P}^{k'}(I')]\in \mathcal{T}_{A_i}, \mbox{ and } \mathcal{P}^{k'}(I')=I\},$$ where $k$ is the generation of $I$, namely, $I\in \mathcal{F}^0_{k}$.
Obviously, the proof of the irreducible case above yields $\dim_H
K_{A_i}(I)\leq \alpha$.

For each $A_i\in \mathcal{E}$, define $\mathcal{F}^0_{A_i}:=\{I\in
\mathcal{F}^0: [I]\in\mathcal{T}_{A_i}\}.$
Then it is easy to verify that
$K=\bigcup_{A_i\in\mathcal{E}}\bigcup_{I\in
\mathcal{F}_{A_i}^0}K_{A_i}(I).$ Hence, it follows from the
countable stability of the Hausdorff dimension (see \cite{Fal}) that
$\dim_{H}K\leq \alpha$, which is the required upper bound.

We have proved that $\mathcal{H}^\alpha(K)>0$ and $\dim_H K=\alpha$.
This imply that $\mathcal{H}^\alpha(K)<\infty$ since $K$ is a
self-similar set (see \cite{Fal}). The proof is completed. $\Box$

We conclude this section with the following examples.

\textbf{Example 2.4.} If $\{S_j\}_{j=1}^m$ satisfies the OSC (with
the open set  (0,1)), then it satisfies the GFTC. The classical
dimension result for $K$ is covered by the result of Theorem 2.3.
See $\cite{Lau}$ for general proof where for each $k\geq 0$,
$\mathcal{M}_k$ is choosen as $\Sigma_k$.

 The following example taken from $\cite{Lau,Wang}$ is an IFS of contractive similitudes whose contraction ratios are not exponentially commensurable. As pointed in $\cite{Lau}$, it satisfies the GFTC, but not  the FTC.

\textbf{Example 2.5.} Let $\{S_j\}_{j=1}^3$ be an IFS on
$\mathbb{R}$ as follows.
$$
S_1(x)=\rho x, \quad S_2(x)=rx+\rho(1-r), \quad S_3(x)=rx+(1-r),
$$
where $0<\rho<1$, $0<r<1$, and $\rho+2r-\rho r\leq 1$. Then
$\{S_j\}_{j=1}^3$ satisfies the GTFC with respect to  $(0,1)$. The
dimension $\alpha$ is the unique solution of the equation
$$\rho^\alpha+2 r^\alpha-(\rho r)^\alpha=1.$$

\setlength{\unitlength}{0.5cm}
\begin{figure}[htbp]
\begin{center}
\includegraphics[width=14cm]{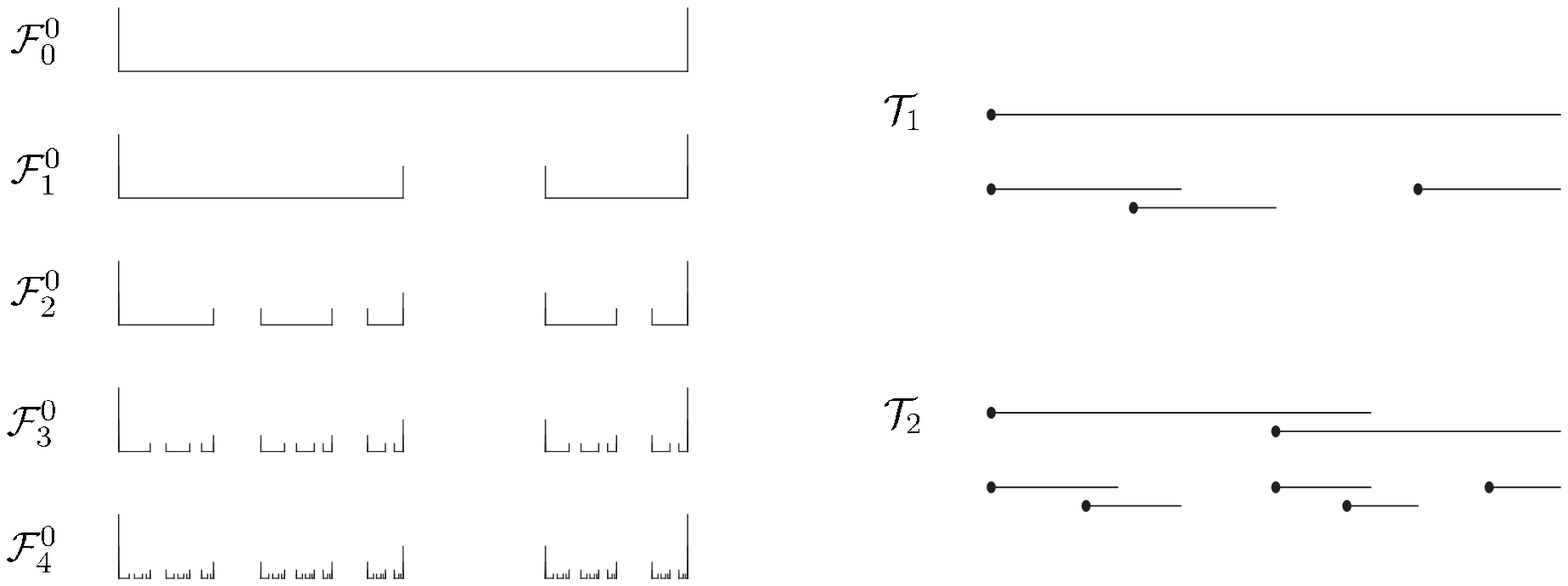}\\
\end{center}\emph{ Figure 1.  The
first five levels of islands and   the distinct overlap types in
Example 2.5 with parameters $\rho=1/3$ and $r=1/4$.}
\end{figure}

\emph{Proof.} For each $k\geq 0$, let $\mathcal{M}_k=\Sigma_k$.
Denote $I_{1}=[0,1]$ and $I_2=S_1([0,1])\cup S_2([0,1])$. It is easy
to verify that $[I_1]$ and $[I_2]$, denoted respectively by
$\mathcal{T}_1$ and $\mathcal{T}_2$, are the total distinct overlap
types. Hence $\mathcal{F}^0/\sim=\{\mathcal{T}_1,\mathcal{T}_2\}$.
The weighted incidence matrix is
$$A_\alpha=\left[\begin{array}{cc}
              r^\alpha & (\rho+r-\rho r)^\alpha  \\
                \frac{r^{2\alpha}}{(\rho+r-\rho r)^\alpha} & \rho^\alpha+r^\alpha \\
                                       \end{array}\right].$$
                                       Setting the spectral radius of $A_\alpha$ equal to $1$ yields the desired result. See Figure 1. $\Box$

In this Example, there exist touching islands if and only if the
contractive ratios $\rho$ and $r$ satisfy the equality $\rho+2r-\rho
r=1$. In this case, the invariant set $K=[0,1]$.

The following example is taken from $\cite{Ngai}$ which is also
satisfying the FTC.

 \textbf{Example
2.6.} Let $\{S_j\}_{j=1}^3$ be an IFS on $\mathbb{R}$ as follows.
$$S_1(x)=\frac{1}{3} x, \quad S_2(x)=\frac{1}{9}x+\frac{8}{27}, \quad S_3(x)=\frac{1}{3}x+\frac{2}{3}.$$
 Then $\{S_j\}_{j=1}^3$ satisfies the GTFC with respect to $(0,1)$. The dimensions
 $\alpha$ ($\approx 0.7369$) is the logarithmic ratio of
the largest root of the polynomial equation
$$x^3-6 x^2+5 x-1=0$$ to $9$.

\setlength{\unitlength}{0.5cm}
\begin{figure}[htbp]
\begin{center}
\includegraphics[width=14cm]{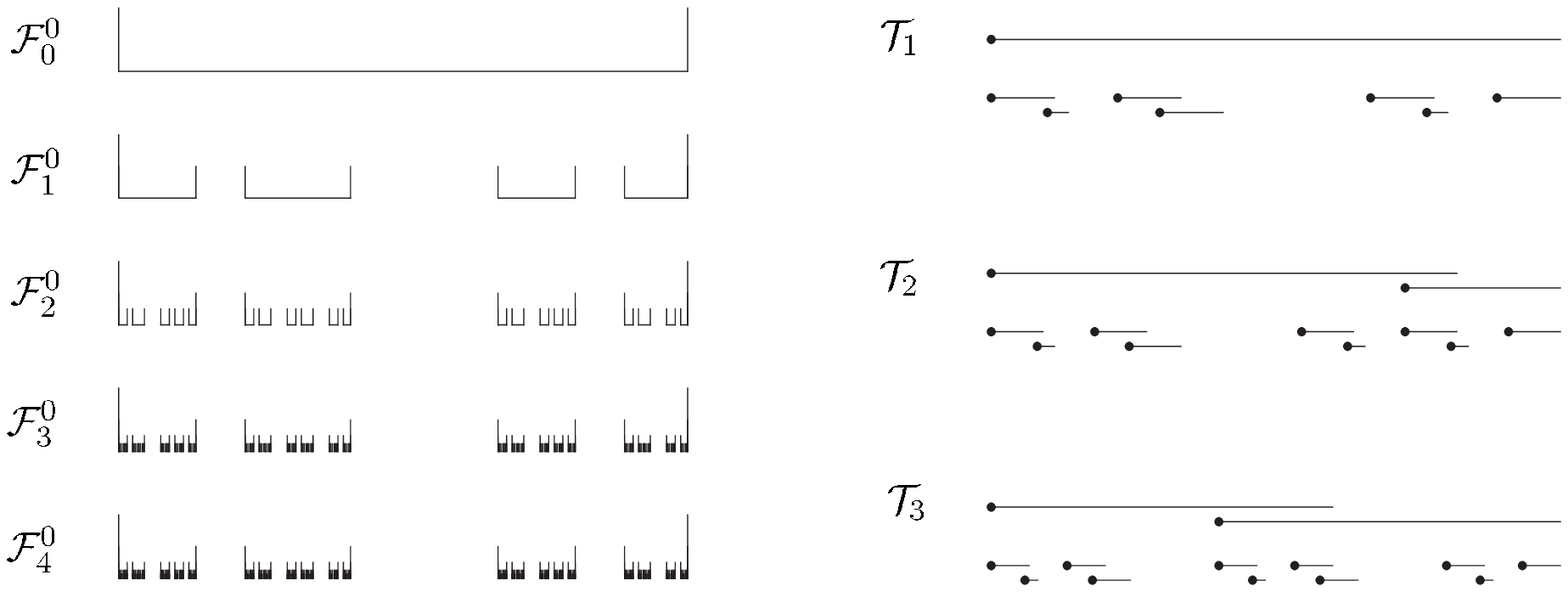}\\
\end{center}\emph{ Figure 2.  The
first five levels of islands and   the distinct overlap types  in
Example 2.6.}
\end{figure}

\emph{Proof.} For each $k\geq 0$, let $\mathcal{M}_k=\Lambda_k$.
Denote $I_1=[0,1]$, $I_2=S_{11}([0,1])\cup S_{12}([0,1])$ and
$I_3=S_{13}([0,1])\cup S_{2}([0,1])$. It is easy to verified that
 $\mathcal{T}_1=[I_1],\mathcal{T}_2=[I_{2}]$ and $\mathcal{T}_3=[I_3]$ are
 the total
overlap types, with a weighted incidence matrix
$$A_\alpha=\left[\begin{array}{ccc}
               \frac{1}{9^\alpha} & \frac{2\cdot 11^\alpha}{81^\alpha} & \frac{5^\alpha}{27^\alpha} \\
             \frac{1}{11^\alpha} & \frac{3}{9^\alpha} & \frac{5^\alpha}{33^\alpha} \\
             \frac{1}{15^\alpha} & \frac{3\cdot 11^\alpha}{135^\alpha}  &\frac{2}{9^\alpha}
                                       \end{array}\right].$$
Setting the spectral radius of $A_\alpha$ equal to $1$ yields the
desired result. See Figure 2. $\Box$

The following is a non-trivial example which allows touching islands.

\textbf{Example 2.7.} Let $\{S_j\}_{j=1}^4$ be an IFS on
$\mathbb{R}$ as follows.
$$S_1(x)=\frac{1}{4} x, \quad S_2(x)=\frac{1}{4}x+\frac{1}{4}, \quad S_3(x)=\frac{1}{4}x+\frac{3}{8}, \quad S_4(x)=\frac{1}{4}x+\frac{3}{4}.$$
 Then $\{S_j\}_{j=1}^4$ satisfies the GTFC with respect to $(0,1)$.
 The dimension is equal to $\alpha={\log_4({5+\sqrt{5}})}-\frac{1}{2}\approx 0.9276.$

\setlength{\unitlength}{0.5cm}
\begin{figure}[htbp]
\begin{center}
\includegraphics[width=14cm]{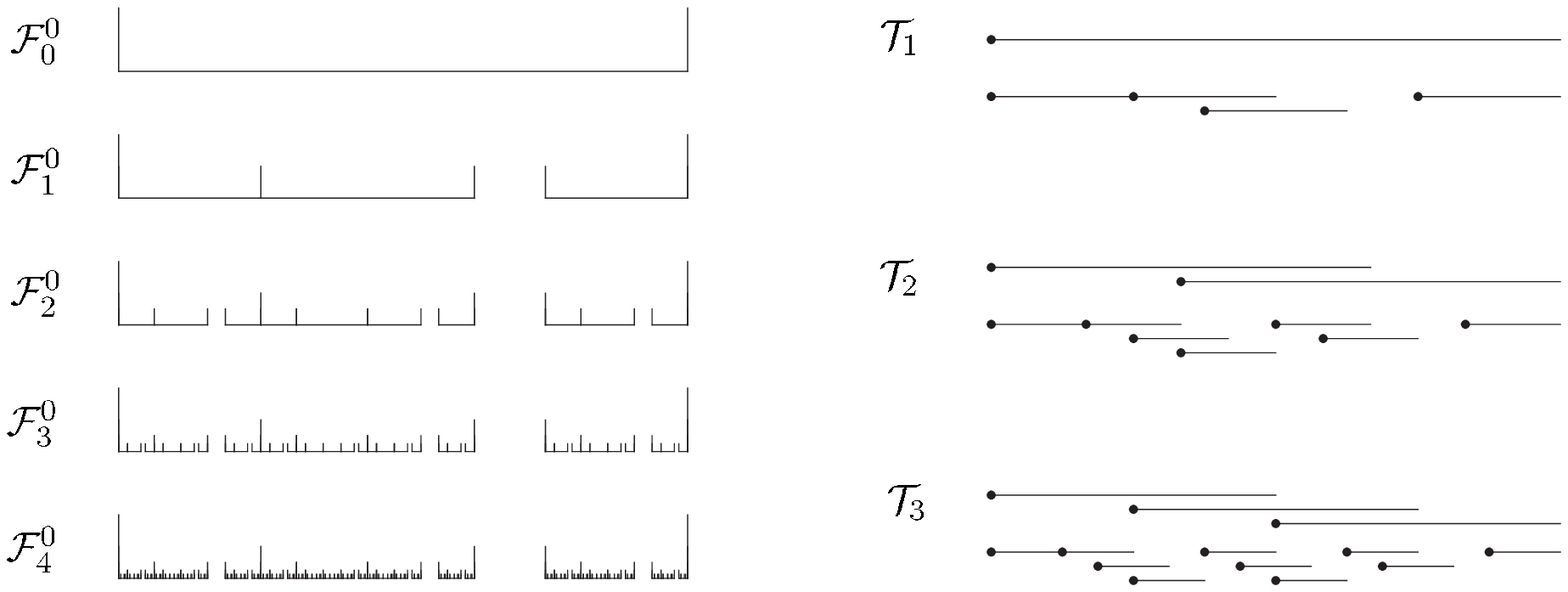}\\
\end{center}\emph{ Figure 3.  The
first five levels of islands and  the distinct overlap types in
Example 2.7.}
\end{figure}

\emph{ Proof. } For each $k\geq 0$, let $\mathcal{M}_k=\Sigma_k$.
Denote $I_1=[0,1]$, $I_2=S_{2}([0,1])\cup S_{3}([0,1])$ and
$I_3=S_{22}([0,1])\cup S_{23}([0,1])\cup S_{31}([0,1])$. It is
straightforward to verify that
$\mathcal{T}_1=[I_1],\mathcal{T}_2=[I_2]$ and $\mathcal{T}_3=[I_3]$
are the total overlap types, with a weighted incidence matrix
$$A_\alpha=\left[\begin{array}{ccc}
               \frac{2}{4^\alpha} & \frac{3^\alpha}{8^\alpha} & 0 \\
             \frac{2}{6^\alpha} & \frac{1}{4^\alpha} & \frac{1}{3^\alpha} \\
             \frac{2}{8^\alpha} & \frac{3^\alpha}{16^\alpha}  &\frac{2}{4^\alpha}
                                       \end{array}\right].$$ Setting the spectral radius of $A_\alpha$ equal to $1$ yields the desired result. See
Figure 3. $\Box$

\section{Density Theorems}

Let $K$ be a linear Cantor set satisfying the GFTC with respect to
the invariant set $\mathcal{O}=(0,1)$.   Let  $\alpha$ be the number
such that the spectral radius of the weighted incidence matrix
$A_\alpha$ is equal to $1$. Then $\alpha$ is the dimension of $K$ by
Theorem 2.3. We assume that the matrix $A_\alpha$ is irreducible
throughout this section and the following ones. It is easy to verify
that all the examples listed in Section 2 satisfy this assumption.
In this section, we analyze the local behavior of the Hausdorff
measure and the packing measure of  $K$.

 For a matrix $A$,
we use $r(A)$ to denote the spectral radius of $A$. The following
basic algebraic lemma is needed.

\textbf{Lemma 3.1.}\emph{ Let $A$ be a $q\times q$ non-negative irreducible matrix with $q\geq 2$, $A'$ be the $(q-1)\times (q-1)$ sub-matrix  at the right-bottom  corner of $A$. Then $r(A')<r(A)$.}

\textit{Proof.} Since $A$ is irreducible, there exists at least one positive number in the set $\{A(1,j): 1\leq j\leq q\}\cup\{A(i,1): 1\leq i\leq q\}$. We define a $q\times q$ matrix $B$ in which $B(1,j)=A(1,j)/2$ for $1\leq j\leq q$, $B(i,1)=A(i,1)/2 $ for $1\leq i\leq q$ and $B(i,j)=A(i,j)$ for $2\leq i,j\leq q$. Obviously, $B$ is also a non-negative irreducible matrix and $B<A$, i.e., their exists at least one coordinate $(i,j)$ such that $B(i,j)<A(i,j)$. Hence from the well-known Perron-Frobenius Theorem, it yields that $r(B)<r(A)$. Then we get the desired result since $r(A')\leq r(B)$. $\Box$

\textbf{Lemma 3.2.} S\emph{uppose $A_\alpha$ is irreducible, then $K$ can be decomposed into an union of a set $K_a$ with a graph directed construction and an attractor $K_b$ of a countable infinite IFS. Moreover, $$\dim_{H}K_a=\dim_P K_a<\dim_H K_b=\dim_P K_b=\dim_H K=\dim_P K,$$  and
 $\mathcal{H}^\alpha(K)=\mathcal{H}^\alpha(K_b), \mathcal{P}^\alpha(K)=\mathcal{P}^\alpha(K_b).$}

\textit{Proof.} For each $k\geq 1$, let $\mathcal{F}'^0_k$ be a subset contained in  $\mathcal{F}^0_k$ as $$\mathcal{F}'^0_k=\{I\in \mathcal{F}^0_k: [I]\neq \mathcal{T}_1\}.$$

Let $K_a=\bigcap_{k\geq 1}\bigcup\{I: I\in\mathcal{F}'^0_k\}$. Then
it is easy to find that $K_a$ is a set with a graph directed
construction. In fact, there are $q-1$ vertex sets in the
construction of $K_a$ whose weighted incidence matrix is the
$(q-1)\times(q-1)$ sub-matrix $A'_\alpha$ at the right-bottom corner
of $A_\alpha$ for each $\alpha\geq 0$. Moreover, by a similar proof
of Theorem 2.3, the Hausdorff dimension of $K_a$ is the unique
$\alpha$ such that the spectral radius of  $A'_\alpha$ is equal to
$1$ (or see a direct result in $\cite{Mau}$). It yields from Lemma
3.1, $\dim_H K_a<\dim_H K$.

Let  $I\in \mathcal{F}^0_k$ with $k\geq 1$. If none of the elements
in $\{\mathcal{P}(I),\cdots,\mathcal{P}^{k-1}(I)\}$ is of
$\mathcal{T}_1$ type, we call $I$ a $\mathcal{T}_1$ \emph{ type
utmost island}. For each $\mathcal{T}_1$ type utmost island $I$, the
vertex set $V(I)$ of  $I$ consists of exactly one vertex, i.e.,
$\sharp V(I)=1$. Denote the contractive similitude of the unique
element in $V(I)$ as $S_I$. Then there exists a countable infinite
IFS (see \cite{Mau1} for further properties of infinite IFS) of
contractive similitudes
$$\mathcal{S}:=\{S_I: I\in \mathcal{F}^0\setminus \mathcal{F}^0_0, I \mbox{ is a } \mathcal{T}_1 \mbox{ type utmost island}\}.$$
Denote by $K_b$ the attractor of $\mathcal{S}$. From the
construction of $K_a$ and $K_b$, one can easily observe that
$K=K_a\cup K_b$. Hence, it follows from the stability property of
the Hausdorff dimension (see \cite{Fal}), $\dim_H K_b=\dim_H K$. The
remaining is obvious. $\Box$

It is worth while to point out that in the above proof we could
replace $\mathcal{T}_1$  by any other overlap type. With appropriate
modifications, we can also prove Lemma 3.2 in a similar way. We will
not go into the details here.

In order to prove Theorem 1.1 and Theorem 1.2, we need a detailed
analysis of the attractor $K_b$ of a countable infinite IFS
$\mathcal{S}$ described in Lemma 3.2. Now we introduce some
notations for convenience. Denote the list of countable contractive
similitudes in $\mathcal{S}$ as
$$\mathcal{S}=\{S'_1,\cdots,S'_j,\cdots\}$$ and $r'_j$  the contractive ratio of $S'_j$. Then for each $k\geq 1$ we will write $S'_\textbf{i}=S'_{i_1}\circ\cdots\circ S'_{i_k}$ and $r'_\textbf{i}=r'_{i_1}\cdots r'_{i_k}$ for all indices $\textbf{i}=i_1\cdots i_k$ with entries $i_j\in \mathbb{N}$. Also, for every such indice $\textbf{i}=i_1\cdots i_k$ we will write $|\textbf{i}|=k$ for the length of $\textbf{i}$.
From the construction of $K_b$, we know the intervals $S'_i([0,1])$ and $S'_j([0,1])$ are separate for all $i\neq j\in \mathbb{N}$. Hence
$$\mathcal{H}^\alpha(K_b)=\sum_{i\in \mathbb{N}}\mathcal{H}^\alpha(S'_i K_b)=\sum_{i\in \mathbb{N}}r'^\alpha_i \mathcal{H}^\alpha(K_b).$$ Combing the above formula with the result in Lemma 3.2 and the fact that $0<\mathcal{H}^\alpha(K)<\infty$ from Theorem 2.3, we have
$\sum_{i\in \mathbb{N}}r'^\alpha_i=1.$ Moreover,
$$\sum_{|\textbf{i}|=k}r'^\alpha_\textbf{i}=\left(\sum_{i\in
\mathbb{N}}r'^\alpha_i\right)^k=1$$ for any integer $k\geq 1$. A
useful measure $\lambda$, which is called $\alpha$-dimensional
normalized Hausdorff measure, defined on $K_b$ by
$\lambda({S'_\textbf{i}(K_b)})=r'^\alpha_\textbf{i}$, then extended
to Borel subsets of $K$, will be used. This is a probability measure
which scales on $S'_\textbf{i}(K_b)$, hence
$\mathcal{H}^\alpha|_{K_b}=\mathcal{H}^\alpha|_{K}=\mathcal{H}^\alpha(K)\lambda$.
Obviously, we also have
$\mathcal{P}^\alpha|_{K_b}=\mathcal{P}^\alpha|_{K}=\mathcal{P}^\alpha(K)\lambda$.

We should point out that  $\lambda$ is equal to the natural
probability measure $\mu$ defined in the proof of Theorem 2.3. This
could be easily verified by showing $\lambda$ and $\mu$ are equal on
cylinder sets in $K_b$, i.e.,
$$\mu(S'_\textbf{i}(K_b))=|S'_\textbf{i}([0,1])|^\alpha
a_1=r'^\alpha_\textbf{i}=\lambda(S'_\textbf{i}(K_b)),$$ where
$a_1=1$ is the first element in the normalized $1$-eigenvector of
$A_\alpha$. Hence for each island $I\in \mathcal{F}^0$, $\lambda(I)$
can be calculated by using the parameters of the IFS of $K$.

For any interval $J\subset [0,1]$ we define the \emph{density} of
$J$ with respect to $\lambda$ as
$d(J)=\frac{\lambda(J)}{|J|^\alpha}.$

Theorem 1.1 and Theorem 1.2 are used to obtain density results for
Hausdorff and packing measures of $K$, namely, Corollary 1.3 and
Corollary 1.4. Now we turn to the proofs.

\subsection{Proof of Theorem 1.1 and Corollary 1.3}

The formula in Theorem 1.1 is analogous to that of a general
self-similar set in $\mathbb{R}^d$ which satisfies the OSC (see
\cite{Ols1}). Using Lemma 3.2, it is now possible to adapt the
techniques there to prove Theorem 1.1. We now show two more lemmas
concerning the measure $\mathcal{H}^\alpha(K_b)$ of the attractor
$K_b$. Theorem 1.1 will be a direct corollary of them.

\textbf{Lemma 3.3.} \emph{Let $K_b$ be the attractor described in
Lemma 3.2, $J\subset [0,1]$ be an interval and $k$ be a positive
integer, then
$\mathcal{H}^\alpha(K_b\cap\bigcup_{|\textbf{i}|=k}S'_\textbf{i}J)\geq
\mathcal{H}^\alpha(K_b\cap J).$}

\textit{Proof.} Since $S'_\textbf{j}K_b\subset\bigcup_{|\textbf{i}|=k}S'_\textbf{i}K_b=K_b$ for all $\textbf{j}$ with $|\textbf{j}|=k$, we have
 $K_b\cap\bigcup_{|\textbf{i}|=k}S'_\textbf{i}J\supset \bigcup_{|\textbf{i}|=k}S'_\textbf{i}(K_b\cap J).$ It follows that
 \begin{eqnarray*}
 \mathcal{H}^\alpha(K_b\cap \bigcup_{|\textbf{i}|=k}S'_\textbf{i} J)&\geq& \mathcal{H}^\alpha\left(\bigcup_{|\textbf{i}|=k}S'_\textbf{i}(K_b\cap J)\right)=\sum_{|\textbf{i}|=k}\mathcal{H}^\alpha\left(S'_\textbf{i}(K_b\cap J)\right)\\&=&\sum_{|\textbf{i}|=k}r'^\alpha_\textbf{i}\mathcal{H}^\alpha(K_b\cap J)=\mathcal{H}^\alpha(K_b\cap J). \quad \Box
 \end{eqnarray*}

The following lemma is a revised version of Theorem 1.1 with $K$
replaced by $K_b$.

\textbf{Lemma 3.4.} \emph{The attractor $K_b$ described in Lemma 3.2
satisfies
\begin{equation}\label{3}
\mathcal{H}^\alpha(K_b\cap J)\leq |J|^\alpha
\end{equation} for all intervals $J\subset [0,1]$.}

\textit{Proof.} In order to reach a contradiction, we assume that
$(\ref{3})$ is not satisfied, i.e., there exists a non-empty
interval $J\subset [0,1]$, such that $\mathcal{H}^\alpha(K_b\cap J)>
|J|^\alpha.$ It follows from this we can find $0<\kappa<1$ with
\begin{equation}\label{4}
(1-\kappa)\mathcal{H}^\alpha(K_b\cap J)>|J|^\alpha.
\end{equation}

Next, fix $\delta>0$ and choose a positive integer $k$ such that
$|S'_\textbf{i}J|\leq \delta$ for all $\textbf{i}$ with
$|\textbf{i}|=k$. Let $\eta=\frac{1}{2}\kappa
\mathcal{H}^\alpha(K_b\cap\bigcup_{|\textbf{i}|=k}S'_\textbf{i}J).$
It follows from Lemma 3.3 and $(\ref{4})$,
\begin{equation}\label{6}
\mathcal{H}^\alpha(K_b\cap\bigcup_{|\textbf{i}|=k}S'_\textbf{i}J)\geq \mathcal{H}^\alpha(K_b\cap J)\geq \frac{|J|^\alpha}{1-\kappa}>0,
\end{equation}  which yields $\eta>0$.

Since $\eta>0$, we can choose a covering $\{J_i\}_i$ of $K_b\setminus \bigcup_{|\textbf{i}|=k}S'_\textbf{i}J$  with $|J_i|\leq \delta$ such that

\begin{eqnarray}\label{5}
\sum_i|J_i|^\alpha&\leq& \mathcal{H}_\delta^\alpha(K_b\setminus
\bigcup_{|\textbf{i}|=k}S'_\textbf{i}J)+\eta\leq
\mathcal{H}^\alpha(K_b\setminus
\bigcup_{|\textbf{i}|=k}S'_\textbf{i}J)+\eta
\end{eqnarray}

The family $\{S'_\textbf{i}J\}_{|\textbf{i}|=k}\cup \{J_i\}_i$ is
clearly a $\delta$-covering of $K_b$. We therefore conclude from
$(\ref{4}), (\ref{5})$ and Lemma 3.3 that
\begin{eqnarray*}
\mathcal{H}^\alpha_\delta(K_b)&\leq& \sum_{|\textbf{i}|=k}|S'_\textbf{i}J|^\alpha+\sum_i |J_i|^\alpha\leq \sum_{|\textbf{i}|=k}r'^\alpha_\textbf{i} |J|^\alpha+\mathcal{H}^\alpha(K_b\setminus \bigcup_{|\textbf{i}|=k}S'_\textbf{i}J)+\eta\\
&=& |J|^\alpha+\mathcal{H}^\alpha(K_b\setminus\bigcup_{|\textbf{i}|=k}S'_\textbf{i}J)+\eta\\&\leq& (1-\kappa)\mathcal{H}^\alpha(K_b\cap\bigcup_{|\textbf{i}|=k}S'_\textbf{i}J)+\mathcal{H}^\alpha(K_b\setminus\bigcup_{|\textbf{i}|=k}S'_\textbf{i}J)+\eta\\
&\leq& \mathcal{H}^\alpha(K_b)-\kappa \mathcal{H}^\alpha(K_b\cap
\bigcup_{|\textbf{i}|=k}S'_\textbf{i}J)+\eta\\&=&
\mathcal{H}^\alpha(K_b)-\eta\leq
\mathcal{H}^\alpha(K_b)-\frac{1}{2}\kappa\mathcal{H}^\alpha(K_b\cap
J).
\end{eqnarray*}

Finally, letting $\delta\rightarrow 0$ gives
$\mathcal{H}^\alpha(K_b)\leq
\mathcal{H}^\alpha(K_b)-\frac{1}{2}\kappa\mathcal{H}^\alpha(K_b\cap
J).$ Then the fact that $0<\mathcal{H}^\alpha(K_b)<\infty$ by Lemma
3.2 and $(1/2)\kappa\mathcal{H}^\alpha(K_b\cap J)>0$ by $(\ref{6})$
provides the desired contradiction. $\Box$

\textit{Proof of Theorem 1.1.} Indeed, it follows from Lemma 3.2 and
Lemma 3.4. $\Box$

 Putting $J$ equal to $[0,1]$ in Theorem 1.1 gives the upper bound for $\mathcal{H}^\alpha(K)$, namely, $\mathcal{H}^\alpha(K)\leq 1$.
  It is a natural generalization of the same result in the OSC case, since in that case $K$ can be covered by its iterated images under the IFS. Obviously it is not always true that we have the equality. See Falconer $\cite{Fal}$ for some examples with $\mathcal{H}^\alpha(K)=1$ satisfying the OSC. An natural question is arisen: \emph{is there any non-trivial example with $\mathcal{H}^\alpha(K)=1$ which satisfies only the GFTC.}

 For a given measure $\nu$ on $\mathbb{R}$ and $s>0$, the \emph{upper $s$-dimensional convex density} of $\nu$ at $x$ is defined by
 $$\mathcal{D}^{*s}(\nu,x)=\lim_{r\rightarrow 0}\sup\{\frac{\nu(J)}{|J|^s}: J \mbox{ is an interval and } 0<|J|\leq r, x\in
 J\}.$$
 The \emph{lower $s$-dimensional convex density} $\mathcal{D}_*^s(\nu,x)$ is defined similarly by taking the lower limit.
We have the following result that if $E\subset \mathbb{R}$ and $s>0$
with $\mathcal{H}^s(E)<\infty$, then
\begin{equation}\label{7}
\mathcal{D}^{*s}(\mathcal{H}^s|_E,x)=1 \mbox{ for } \mathcal{H}^s\mbox{-a.e. } x\in E.
\end{equation}
The reader is referred to $\cite{Fal1}$ for a proof of $(\ref{7})$.

\textit{Proof of Corollary 1.3.}
From $(\ref{7})$, we can pick a point $x\in K\cap (0,1)$ such that $\mathcal{D}^{*\alpha}(\mathcal{H}^\alpha|_K,x)=1$. By the definition of $\mathcal{D}^{*\alpha}(\mathcal{H}^\alpha|_K,x)$, there exists a positive sequence $\{\delta_n\}_n$ with $\delta_n<\min\{x,1-x\}$ and $\delta_n\rightarrow 0$ as $n\rightarrow \infty$ such that
$$\sup_{0<|J|\leq\delta_n}\frac{\mathcal{H}^\alpha(K\cap J)}{|J|^\alpha}-\frac{1}{n}\leq 1\leq \sup_{0<|J|\leq\delta_n}\frac{\mathcal{H}^\alpha(K\cap J)}{|J|^\alpha}+\frac{1}{n}.$$
Hence there exists an interval $J_n$ with $0<|J_n|\leq \delta_n$ for each $n$ such that
$$\sup_{0<|J|\leq\delta_n}\frac{\mathcal{H}^\alpha(K\cap J)}{|J|^\alpha}\leq \frac{\mathcal{H}^\alpha(K\cap J_n)}{|J_n|^\alpha}+\frac{1}{n}.$$
Thus $\frac{\mathcal{H}^\alpha(K\cap
J_n)}{|J_n|^\alpha}-\frac{1}{n}\leq
1\leq\frac{\mathcal{H}^\alpha(K\cap
J_n)}{|J_n|^\alpha}+\frac{2}{n},$ which yields that
$\frac{\mathcal{H}^\alpha(K\cap J_n)}{|J_n|^\alpha}\rightarrow 1$ as
$n\rightarrow \infty$. Moreover, by Theorem 1.1, for each interval
$J_n$ we have ${\mathcal{H}^\alpha(K\cap J_n)}/{|J_n|^\alpha}\leq
1$. Hence $\sup\{\frac{\mathcal{H}^\alpha(K\cap J)}{|J|^\alpha}:
J\subset [0,1]\}=1.$ Since
$\lambda={\mathcal{H}^\alpha|_K}/{\mathcal{H}^\alpha(K)}$,
$(\ref{8})$ follows immediately from the above equation. $\Box$

 With
suitable modifications if necessary, we may generalize Theorem 1.1
and Corollary 1.3 in an similar way to general self-similar sets
satisfying the GFTC with irreducible weighted incidence matrix.
However, in order to match the main goal of this paper and to avoid
additional technical details, we will not pursue this here.

\subsection{Proof of Theorem 1.2 and Corollary 1.4}

 In a manner dual to the Hausdorff measure case, the packing measure result is also analogous to that of a general self-similar set in $\mathbb{R}^d$ which satisfies the SSC or the OSC (see \cite{Ols1,Qiu}). Hence it is also possible to adapt the techniques there to prove Theorem 1.2 by using Lemma 3.2.
We shall need the following two lemmas  concerning
$\mathcal{P}^\alpha(K_b)$.

\textbf{Lemma 3.5.} \textit{Let $K_b$ be the attractor described in
Lemma 3.2. Let $J\subset[0,1]$ be an  interval centered in $K_b$ and
$k$ a positive integer. Then $ \mathcal{P}^\alpha(K_b\cap
\bigcup_{|\textbf{i}|=k}S'_{\textbf{i}}J)=\mathcal{P}^\alpha(K_b\cap
J)>0. $}

\textit{Proof.} We write $J^\circ$ as the interior of the interval $J$. First, we prove that
\begin{equation}\label{10}
K_b\cap\bigcup_{|\textbf{i}|=k}S'_\textbf{i}J^\circ=\bigcup_{|\textbf{i}|=k}S'_\textbf{i}(K_b\cap J^\circ).
\end{equation}

Fix a point $y\in
K_b\cap\bigcup_{|\textbf{i}|=k}S'_\textbf{i}J^\circ$. Since
$J^\circ\subset \mathcal{O}=(0,1)$, there exists an index
$\textbf{u}$ with the length $k$ such that $y\in
S'_\textbf{u}J^\circ\subset S'_\textbf{u}\mathcal{O}$. We also have
$y\in K_b=\bigcup_{|\textbf{i}|=k}S'_\textbf{i}K_b$ and we therefore
find an index $\textbf{v}$ with the length $k$ such that $y\in
S'_\textbf{v}K_b\subset S'_\textbf{v}\overline{\mathcal{O}}$. Thus
$y\in S'_\textbf{u}\mathcal{O}\cap
S'_\textbf{v}\overline{\mathcal{O}}$, and therefore
$\textbf{u}=\textbf{v}.$ Hence $y\in S'_\textbf{u}J^\circ\cap
S'_\textbf{u}K_b=S'_\textbf{u}(K_b\cap
J^\circ)\subset\bigcup_{|\textbf{i}|=k}S'_\textbf{i}(K_b\cap
J^\circ)$ which yields that
$K_b\cap\bigcup_{|\textbf{i}|=k}S'_\textbf{i}J^\circ\subset\bigcup_{|\textbf{i}|=k}S'_\textbf{i}(K_b\cap
J^\circ).$  The other direction is obvious.

It follows from (\ref{10}) that
\begin{eqnarray*}
\mathcal{P}^\alpha(K_b\cap\bigcup_{|\textbf{i}|=k}S'_\textbf{i}J)
&=&\mathcal{P}^\alpha\left(\bigcup_{|\textbf{i}|=k}S'_\textbf{i}(K_b\cap
J^\circ)\right)=\sum_{|\textbf{i}|=k}\mathcal{P}^\alpha\left(S'_\textbf{i}(K_b\cap
J^\circ)\right)\\&=&\sum_{|\textbf{i}|=k}r'^\alpha_\textbf{i}
\mathcal{P}^\alpha(K_b\cap J^\circ)=\mathcal{P}^\alpha(K_b\cap J).
\end{eqnarray*}

Moreover, since $J$ has its center in $K_b$, we deduce that
$\mathcal{P}^\alpha(K_b\cap J)>0$. This completes the proof of Lemma
3.5. $\Box$

\textbf{Lemma 3.6.} The attractor $K_b$ described in Lemma 3.2
satisfies
\begin{equation}\label{99}
\mathcal{P}^\alpha(K_b\cap J)\geq |J|^\alpha
\end{equation} for all intervals $J\subset [0,1]$ centered in $K_b$.

\textit{Proof.} In order to reach a contradiction, we assume that
$(\ref{99})$ is not satisfied, i.e., there exists an interval
$J:=[c,d]\subset [0,1]$ centered in $K_b$, such that
$\mathcal{P}^\alpha(K_b\cap J)<|J|^\alpha.$ Thus we can find a
number $0<\kappa<1$ with
\begin{equation}\label{11}
(1+\kappa)\mathcal{P}^\alpha(K_b\cap J)<|J|^\alpha.
\end{equation}

Next, fix $\delta>0$ and choose a positive integer $k$ such that
$|S'_\textbf{i}J|\leq \delta$ for all $\textbf{i}$ with
$|\textbf{i}|=k$. Let
$\eta=\frac{1}{2}\kappa\mathcal{P}^\alpha(K_b\cap J).$ It follows
from Lemma 3.5, $\eta>0.$

For a positive integer $n,$ write
$G_n=K_b\setminus\bigcup_{|\textbf{i}|=k}S'_\textbf{i}\left((c-{1}/{n},d+{1}/{n})\right),$
and observe that $G_1\subset G_2\subset\cdots\subset
G_n\subset\cdots$ and $\bigcup_n
G_n=K_b\setminus\bigcup_{|\textbf{i}|=k}S'_\textbf{i}J.$

If $K_b\setminus\bigcup_{|\textbf{i}|=k}S'_\textbf{i}J\neq
\emptyset,$ then there is a positive integer ${n_0}$ with
${1}/{{n_0}}<\delta$ such that $G_{n_0}\neq \emptyset$, and
\begin{equation}\label{12}
\mathcal{P}^\alpha(G_{n_0})\geq \mathcal{P}^\alpha(K_b\setminus\bigcup_{|\textbf{i}|=k}S'_\textbf{i}J)-\frac{\eta}{2}.\end{equation}

We can choose a ${1}/{{n_0}}$-packing $\{J^\circ_i\}_i$ of $G_{n_0}$ such that
\begin{equation}\label{13}
\sum_{i}|J_i|^\alpha\geq
P^\alpha_{\frac{1}{{n_0}}}(G_{n_0})-\frac{\eta}{2}\geq
P^\alpha(G_{n_0})-\frac{\eta}{2}\geq
\mathcal{P}^{\alpha}(G_{n_0})-\frac{\eta}{2}.
\end{equation}

Since $J^\circ\subset \mathcal{O}$, $S'_\textbf{i}(J^\circ)\cap
S'_\textbf{j}(J^\circ)=\emptyset$ for all $\textbf{i}\neq
\textbf{j}$ with lengths $k$. And for each such $\textbf{i}$, since
$({c+d})/{2}\in K_b$, we have
$S'_\textbf{i}\left(({c+d})/{2}\right)\in S'_\textbf{i}K_b\subset
K_b$. Thus the family $\{S'_\textbf{i}(J^\circ)\}_{|\textbf{i}|=k}$
is a $\delta$-packing of
$K_b\cap\bigcup_{|\textbf{i}|=k}S'_\textbf{i}J$ by the fact that
$|S'_\textbf{i}J|\leq\delta$ for all $\textbf{i}$ with
$|\textbf{i}|=k$.

Since $\{J^\circ_i\}_i$ is also  a ${1}/{{n_0}}$-packing of
$G_{n_0}$, we conclude that
$\{S'_\textbf{i}J^\circ\}_{|\textbf{i}|=k}\bigcup\{J_i^\circ\}_i$ is
a $\delta$-packing of $K_b$. Using this we therefore conclude from
$(\ref{11})$, $(\ref{12})$, $(\ref{13})$, and Lemma 3.5 that
\begin{eqnarray*}
P^\alpha_\delta(K_b)&\geq& \sum_{|\textbf{i}|=k}(r'_\textbf{i}|J|)^\alpha+\sum_i(|J_i|)^\alpha\geq \sum_{|\textbf{i}|=k}r'^\alpha_\textbf{i}|J|^\alpha+\mathcal{P}^{\alpha}(G_{n_0})-\frac{\eta}{2}\\&\geq& |J|^\alpha+\mathcal{P}^\alpha(K_b\setminus\bigcup_{|\textbf{i}|=k}S'_\textbf{i}J)-\eta\\
&\geq &(1+\kappa)\mathcal{P}^\alpha(K_b\cap J)+\mathcal{P}^\alpha(K_b\setminus\bigcup_{|\textbf{i}|=k}S'_\textbf{i}J)-\eta\\
&=& (1+\kappa)\mathcal{P}^\alpha(K_b\cap\bigcup_{|\textbf{i}|=k}S'_\textbf{i}J)+\mathcal{P}^\alpha(K_b\setminus\bigcup_{|\textbf{i}|=k}S'_\textbf{i}J)-\eta\\
&=& \mathcal{P}^\alpha(K_b)+\kappa\mathcal{P}^\alpha(K_b\cap\bigcup_{|\textbf{i}|=k}S'_\textbf{i}J)-\eta\\
&=& \mathcal{P}^\alpha(K_b)+\frac{1}{2}\kappa \mathcal{P}^\alpha(K_b\cap J).
\end{eqnarray*}

Finally, let $\delta\rightarrow 0$, we get
\begin{equation}\label{14}
P^\alpha(K_b)\geq \mathcal{P}^\alpha(K_b)+\frac{1}{2}\kappa\mathcal{P}^\alpha(K_b\cap J).
\end{equation}

In $\cite{Fen}$ it is proved that the packing premeasure $P^\alpha$
coincides with the packing measure $\mathcal{P}^\alpha$ for compact
subsets with finite $P^\alpha$-measure. Thus they coincide for
$K_b$, and it follows from $(\ref{14})$ that
$\mathcal{P}^\alpha(K_b)\geq
\mathcal{P}^\alpha(K_b)+\frac{1}{2}\kappa\mathcal{P}^\alpha(K_b\cap
J).$ Since $\mathcal{P}^\alpha(K_b)$ is positive and finite by Lemma
3.2, and $({1}/{2})\kappa\mathcal{P}^\alpha(K_b\cap J)>0$, we get
the contradiction.

On the other hand, if
$K_b\setminus\bigcup_{|\textbf{i}|=k}S'_\textbf{i}J=\emptyset,$
i.e., $K_b\subset\bigcup_{|\textbf{i}|=k}S'_\textbf{i}J$, the
aforementioned string of  inequalities simplifies to
$$
P^\alpha_\delta(K_b)\geq
\sum_{|\textbf{i}|=k}(r'_\textbf{i}|J|)^\alpha=
|J|^\alpha>(1+\kappa)\mathcal{P}^\alpha(K_b\cap J)=
(1+\kappa)\mathcal{P}^\alpha(K_b\cap\bigcup_{|\textbf{i}|=k}S'_\textbf{i}J).
$$
Letting $\delta\rightarrow 0$ and using the fact that $\mathcal{P}^\alpha(K_b)=P^\alpha(K_b)$ gives
$$
\mathcal{P}^\alpha(K_b)=P^\alpha(K_b)\geq
(1+\kappa)\mathcal{P}^\alpha(K_b\cap\bigcup_{|\textbf{i}|=k}S'_\textbf{i}J)=(1+\kappa)\mathcal{P}^\alpha(K_b)>\mathcal{P}^\alpha(K_b).
$$
This provides the desired contradiction. $\Box$

\emph{Proof of Theorem 1.2.} For any $J:=[c,d]\subset [0,1]$
centered in $K$. If $J$ is centered in $K_b$, the result follows
immediately from Lemma 3.2 and Lemma 3.6. If $J$ is not centered in
$K_b$, i.e., $(c+d)/2\in K_a\setminus K_b$, we claim that:

\emph{ Claim:} $K_b\cap(c,d)\neq\emptyset$ and for each positive integer $n$, there exists a point $x_n\in K_b\cap (c,d)$ such that $|x_n-(c+d)/2|\leq {1}/{n}$.

\emph{Proof.} Since $(c+d)/2\in K_a$, for each $n$, there is an island $I_n\subset (c,d)$ containing the point $(c+d)/2$ with length less than ${1}/{n}$ whose overlap type is not $\mathcal{T}_1$. By the irreducible property of the weighted incidence matrix $A_\alpha$, we could find a smaller $\mathcal{T}_1$ type island $\widetilde{I}_n$ contained in $I_n$. By the construction of $K_b$, $K_b\cap \widetilde{I}_n\neq\emptyset$ which yields that $K_b\cap(c,d)\neq\emptyset$. Then fix a point $x_n$ in $K_b\cap \widetilde{I}_n$. Obviously we get $|x_n-(c+d)/2|\leq |I_n|\leq {1}/{n}$. $\Box$

By this claim we define a sequence of intervals $\{J_n\}_n$ contained in $J$ as $$J_n=[x_n-\min\{x_n-c,d-x_n\},x_n+\min\{x_n-c,d-x_n\}].$$
It is not difficult to find that for each $n$ the interval $J_n$ is centered in $K_b$, and that the left endpoint of $J_n$ tends to $c$ and the right endpoint of $J_n$ tends to $d$ as $n\rightarrow \infty$. As showed in the first case, $\mathcal{P}^\alpha(K\cap J_n)\geq |J_n|^\alpha$ for all $n$. Letting $n\rightarrow\infty$, we immediately get $\mathcal{P}^\alpha(K\cap J)\geq |J|^\alpha$. $\Box$

Let $s\geq 0$, for a given measure $\nu$ on $\mathbb{R}$ and $x\in \mathbb{R}$, the \emph{lower $s$-dimensional density} of $\nu$ at $x$ is defined as
$$\Theta_*^s(\nu,x)=\liminf_{r\rightarrow 0}\frac{\nu([x-r,x+r])}{(2r)^s}.$$
The \emph{upper $s$-dimensional density} $\Theta^{*\alpha}(\nu,x)$
is defined similarly by taking the upper limit. We have the
following result. If $E\subset \mathbb{R}$ and $s>0$ with
$0<\mathcal{P}^s(E)<\infty$,
\begin{equation}\label{15}
\Theta_*^s(\mathcal{P}^s|_E,x)=1 \mbox{ for } \mathcal{P}^\alpha\mbox{-a.e.} \quad x\in E.
\end{equation}
See the proof in \cite{Mat1}. Now we prove Corollary 1.4 on the
basis of $(\ref{15})$ and Theorem 1.2.

\textit{Proof of Corollary 1.4.} Let $\mathcal{O}=(0,1)$. Since
$K\cap \mathcal{O}\neq\emptyset,$ we can take a point $y\in K\cap
\mathcal{O}$. Choose $\delta>0$ such that the open interval
$(y-\delta,y+\delta)$ is contained in $\mathcal{O}$. Moreover, since
$y\in K$, $\mathcal{P}^\alpha(K\cap (y-\delta,y+\delta))>0.$ Hence
from $(\ref{15})$, there exists a point $z\in K\cap
(y-\delta,y+\delta)$ with
$\Theta_*^\alpha(\mathcal{P}^\alpha|_K,z)=1.$ By the definition of
$\Theta_*^\alpha(\mathcal{P}^\alpha|_K,z)$, there exists a sequence
$\{r_n\}_n$ with each $r_n<\delta-d(z,y)$ and $r_n\rightarrow 0 $ as
$n\rightarrow\infty,$ such that
$\lim_{n\rightarrow\infty}{\mathcal{P}^\alpha(K\cap
[z-r_n,z+r_n])}/{(2r_n)^\alpha}=1.$ Notice that all  intervals
$[z-r_n,z+r_n]$ are contained in
$(y-\delta,y+\delta)\subset\mathcal{O}$ with center $z\in K$.
However, by Theorem 1.2, for each interval $J\subset [0,1]$ centered
in $K$, we have ${\mathcal{P}^\alpha(K\cap J)}/{|J|^\alpha}\geq 1$.
Hence we get $\inf\{\frac{\mathcal{P}^\alpha(K\cap J)}{|J|^\alpha}:
J \mbox{ is an interval centered in } K \mbox{ with } J\subset
[0,1]\}=1.$ Since
$\lambda={\mathcal{P}^\alpha|_K}/{\mathcal{P}^\alpha(K)}$,
$(\ref{16})$ follows immediately from the above equation. $\Box$

 Similar to the
Hausdorff measure case, with suitable modifications if necessary, we
may generalize Theorem 1.2 and Corollary 1.4 in a similar way to
general self-similar sets in $\mathbb{R}^d$ satisfying the GFTC with
irreducible weighted incidence matrix. Due to the same reason for
the Hausdorff measure case, we will not pursue such generalizations
here.

\section{The maximal and minimal densities}

In this section we deal with the exact computation of the Hausdorff
measure and the packing measure for a special kind of linear Cantor
sets. Let $S_j(x)=\rho_jx+b_j, j=1,\cdots,m$, be a linear
contractive IFS on the line $\mathbb{R}$ satisfying the GFTC with
respect to the open set $\mathcal{O}=(0,1)$ with an irreducible
weighted incidence matrix $A_\alpha$. As before, we write $K$ as its
invariant set and use $\alpha$ to denote the dimension of $K$. As
$\cite{Ayer,Feng}$, we do not allow negative $\rho$'s, namely, we
assume $0<\rho_j<1$ for $j=1,\cdots,m$. Moreover, without loss of
generality and for convenience  we assume the images $S_j([0,1])$
are in increasing order, with $S_1(0)=0$ and $S_m(1)=1$. To avoid
triviality we always assume $m\geq 2$ and $\alpha<1$.

Define $$d_{\max}:=\sup\{d(J): J \mbox{ is an interval with } J\subset [0,1]\}$$ the maximal density of intervals contained in $[0,1]$ and
$$d_{\min}:=\inf\{d(J): J \mbox{ is an interval centered in } K \mbox{ with } J\subset [0,1]\}$$ the minimal density of intervals centered in $K$ and contained in $[0,1]$. Then by Corollary 1.3 and Corollary 1.4, $$\mathcal{H}^\alpha(K)=d^{-1}_{\max}\mbox{ and } \mathcal{P}^\alpha(K)=d^{-1}_{\min}.$$ Hence our main purpose in this section is to determine the constants $d_{\max}$ and $d_{\min}$.

We will frequently use the notation $\mathcal{F}_k$, the
\emph{$k$-th generation field}, which is the finite field generated
by the set $\mathcal{F}^0_k$ for $k\geq 0$.  For each $k\geq 0$, let
$\beta_1^{(k)},\cdots,\beta_{l_k}^{(k)}$ denote the lengths of the
$k$-th generation islands in increasing order where $l_k$ denotes
the cardinality of $\mathcal{F}^0_k$. Denote
$\beta^{(k)}_{\max}=\max\{\beta_1^{(k)},\cdots,\beta_{l_k}^{(k)}\}$
and
$\beta^{(k)}_{\min}=\min\{\beta_1^{(k)},\cdots,\beta_{l_k}^{(k)}\}$.
We write the lengths of the lakes separating the $k$-th generation
islands as $\gamma^{(k)}_1,\cdots,\gamma^{(k)}_{l_k-1}$ in
increasing order with
$\gamma_{\min}^{(k)}=\min\{\gamma^{(k)}_1,\cdots,\gamma^{(k)}_{l_k-1}\}$.
Note that we allow $\gamma_{\min}^{(k)}=0$ in the case that touching
islands exist. Denote by $\gamma_{\min}^{*(k)}$ the minimal length
of non-empty lakes in $\mathcal{F}_k$. The identity
$$\beta^{(k)}_1+\cdots+\beta^{(k)}_{l_k}+\gamma^{(k)}_{1}+\cdots+\gamma^{(k)}_{l_k-1}=1,$$ the positivity of $\beta$'s and the  non-negativity of $\gamma$'s are the only restrictions on these parameters. It should be mentioned here that $l_k\rightarrow \infty$, $\beta_{\max}^{(k)}\rightarrow 0$, $\beta_{\min}^{(k)}\rightarrow 0$, $\gamma_{\min}^{(k)}\rightarrow 0$ and $\gamma_{\min}^{*(k)}\rightarrow 0$ as $k\rightarrow \infty$.

The  measure $\lambda$ will play a key role in this section. It
enables us to compute the density of each interval in
$\mathcal{F}_k$ as an elementary function on parameters of the IFS
of $K$. Moreover, there is an obvious algorithm for finding the
maximal or minimal density of intervals in $\mathcal{F}_k$.

Since $K$ has finite types of islands, there exists a smallest
non-negative integer $k_0$ such that none of the islands in
$\mathcal{F}^0_{k_0+1}$ is of a new overlap type. We have the
following \emph{blow-up principle} for the density of each interval
$J\subset [0,1]$.

\textbf{Lemma 4.1.} \emph{For any interval $J$, there exists another
interval $J'$, not contained in a $(k_0+1)$-th generation island,
with the same density.}

\emph{Proof.} If $J$ is contained in a $(k_0+1)$-th generation
island $I$, then there exists a larger island $I'\in \bigcup_{0\leq
k\leq k_0}\mathcal{F}^0_k$ with the same type as $I$, i.e.,
$[I']=[I]$. Let $\tau$ be the linear function which maps $I$ onto
$I'$, keeping the orientation. Then obviously the image $J':=\tau
(J)$ has the same density as $J$. We iterate this procedure until we
obtain $J'$ not lying in any $(k_0+1)$-th generation island. $\Box$

Thus we only need to consider intervals not contained in
$(k_0+1)$-th generation islands.

In order to get general results for computing
$\mathcal{H}^\alpha(K)$ and $\mathcal{P}^\alpha(K)$, we need the
following two additional technical assumptions. We will show that
these assumptions are as general as to be satisfied by all the
examples illustrated in Section 2.

\textbf{Assumption A.} For each island $I\in \mathcal{F}^0$, for any
two different constitutive intervals
$\overline{\mathcal{O}}_\textbf{v}$ and
$\overline{\mathcal{O}}_\textbf{u}$ of $I$,
$\overline{\mathcal{O}}_\textbf{v}\nsubseteq
\overline{\mathcal{O}}_\textbf{u}$ and
$\overline{\mathcal{O}}_\textbf{u}\nsubseteq
\overline{\mathcal{O}}_\textbf{v}$.

\textbf{Assumption B.} $S_{1}([0,1])\cap K=S_1 K$ and
$S_m([0,1])\cap K=S_m K$.

Before we give some remarks on these assumptions, we should prove
that

\textbf{Proposition 4.2.} For all the examples listed in Section 2,
the above two assumptions hold.

\emph{Proof.} It is trivial for Example 2.4 and 2.7.

For Example 2.5,  we only need to prove $S_{1}([0,1])\cap K=S_1 K$.

By the fact $S_{1}([0,1])\cap S_3 K=\emptyset$, $S_1([0,1])\cap
S_{23}K=\emptyset$ and $S_{13}=S_{21}$, we have
\begin{eqnarray*}
S_{1}([0,1])\cap K&=&\left(S_1([0,1])\cap S_1 K\right)\cup \left(S_1([0,1])\cap S_2K\right)\\&=&\left(S_1([0,1]\cap K)\right)\cup \left(S_1([0,1])\cap S_{13}K\right)\cup (S_1([0,1])\cap S_{2}^2K)\\
&=&(S_1([0,1]\cap K))\cup (S_1([0,1]\cap S_{3}K))\cup (S_1([0,1])\cap S_{2}^2K).
\end{eqnarray*}
For $n\geq 2$, iterate the above proceedure $n-1$ times, we get
\begin{eqnarray*}
S_{1}([0,1])\cap K&=&(S_1([0,1]\cap K))\cup (S_1([0,1]\cap S_{3}K))\cup\\&&\cdots\cup (S_{1}([0,1]\cap S_3^nK))\cup (S_1([0,1])\cap S_{2}^{n+1}K).
\end{eqnarray*}
Hence we get  $S_{1}([0,1])\cap K\subset S_1K$ by the arbitrariness
of $n$. The other direction is obvious.

For Example 2.6, we still only need to prove $S_{1}([0,1])\cap K=S_1
K$.

By the fact $S_{1}([0,1])\cap S_3 K=\emptyset$, $S_1([0,1])\cap
S_{23}K=\emptyset$ and $S_{133}=S_{21}$, we have
\begin{eqnarray*}
S_{1}([0,1])\cap K&=&(S_1([0,1])\cap S_1 K)\cup (S_1([0,1])\cap S_2K)\\&=&(S_1([0,1]\cap K))\cup (S_1([0,1])\cap S_{133}K)\cup (S_1([0,1])\cap S_{2}^2K)\\
&=&(S_1([0,1]\cap K))\cup (S_1([0,1]\cap S_{3}^2K))\cup (S_1([0,1])\cap S_{2}^2K).
\end{eqnarray*}
For $n\geq 2$, iterate the above proceedure $n-1$ times, we get
\begin{eqnarray*}
S_{1}([0,1])\cap K&=&(S_1([0,1]\cap K))\cup (S_1([0,1]\cap S_{3}^2K))\cup\\&&\cdots\cup (S_{1}([0,1]\cap S_3^{2n}K))\cup (S_1([0,1])\cap S_{2}^{n+1}K).
\end{eqnarray*}
Hence we get $S_{1}([0,1])\cap K\subset S_1K$ by the arbitrariness
of $n$. The other direction is obvious. $\Box$

\textbf{Remark 1.} Under Assumption B, for  each positive integer
$k$, we have $S_{1}^k([0,1])\cap K=S_1^k K$ and $S_m^{k}([0,1])\cap
K=S_{m}^{k}K$.

\emph{Proof.} We only need to prove $S_{1}^2([0,1])\cap K=S_{1}^2K$.
Since $S_{1}^2([0,1])\cap K\subset S_{1}([0,1])\cap K=S_1 K$, we
have $S_{1}^2([0,1])\cap K=S_{1}^2([0,1])\cap S_1
K=S_1(S_1([0,1])\cap K)=S_{1}^2K$. The second equality  can be
proved in a similar way. $\Box$

\textbf{Remark 2.} For each island $I=[c,d]\in \mathcal{F}^0$,
notice that $I=\bigcup _{\textbf{v}\in
V(I)}\overline{\mathcal{O}}_\textbf{v}$. By Assumption A, there is a
unique vertex $\textbf{v}_0\in V(I)$ such that
$\overline{\mathcal{O}}_{\textbf{v}_0}$ has the same left endpoint
as $I$, i.e., $S_{\textbf{v}_0}(0)=c$. If $V(I)\setminus
\{\textbf{v}_0\}\neq\emptyset$, let $c'=\min\{S_\textbf{v}(0):
\textbf{v}\in V(I)\setminus\{\textbf{v}_0\}\}$. Obviously $c<c'<d$.
Let $k$ be the smallest positive integer such that
$S_{\textbf{v}_0}S_{\textbf{i}_0}([0,1])\subset [c,c')$ where
$\textbf{i}_0=(1,\cdots,1)$ with $|\textbf{i}_0|=k$ (This must be
done since $S_{\textbf{v}_0}S_{1}^k(0)$ is always equal to c). For
this $\textbf{i}_0$, by Remark 1, we have
$$S_{\textbf{v}_0}S_{\textbf{i}_0}([0,1])\cap K=S_{\textbf{v}_0}S_{\textbf{i}_0} K.$$ This can be verified since
$$
S_{\textbf{v}_0}S_{\textbf{i}_0}([0,1])\cap
K=S_{\textbf{v}_0}S_{\textbf{i}_0}([0,1])\cap S_{\textbf{v}_0}
K=S_{\textbf{v}_0}(S_{\textbf{i}_0}([0,1])\cap
K)=S_{\textbf{v}_0}S_{\textbf{i}_0} K. $$ Otherwise, if
$V(I)\setminus \{\textbf{v}_0\}=\emptyset$, i.e.,
$V(I)=\{\textbf{v}_0\}$, we define $\textbf{i}_0=\emptyset$. Also,
we have $S_{\textbf{v}_0}S_{\textbf{i}_0}([0,1])\cap
K=S_{\textbf{v}_0}S_{\textbf{i}_0} K.$

Similarly,  In an analogous way, there is also a unique vertex
$\textbf{v}_1\in V(I)$ such that
$\overline{\mathcal{O}}_{\textbf{v}_1}$ has the same right endpoint
as $I$, i.e., $S_{\textbf{v}_1}(1)=d$.  By a similar discussion, one
can find that there also exists a smallest non-negative integer $k'$
such that
$$S_{\textbf{v}_1}S_{\textbf{i}_1}([0,1])\cap K=S_{\textbf{v}_1}S_{\textbf{i}_1} K$$ where $\textbf{i}_1=(m,\cdots,m)$ with the length $|\textbf{i}_1|=k'$.

\textbf{Remark 3.} We should mention that in Remark 2 the ratios
$|S_{\textbf{v}_0}([0,1])|/|I|,$ $|S_{\textbf{v}_1}([0,1])|/|I|$ and
the indices $\textbf{i}_0, \textbf{i}_1$ are dependent merely on the
overlap type of $I$ and are independent on the choice of $I$. In
order to emphasize the relation between
$\textbf{v}_0,\textbf{v}_1,\textbf{i}_0,\textbf{i}_1$ and $I$, we
replace $\textbf{v}_0,\textbf{v}_1,\textbf{i}_0,\textbf{i}_1$ by
$\textbf{v}_0(I),\textbf{v}_1(I),\textbf{i}_0(I),\textbf{i}_1(I)$
respectively. Since the number of the overlap types is finite, we
could define two positive numbers $\eta_1\leq 1$ and $\eta_2\leq 1$
as follows.
$$\eta_1:=\min_{I\in\mathcal{F}^0}\{\frac{\rho_{\textbf{v}_0(I)}}{|I|}, \frac{\rho_{\textbf{v}_1(I)}}{|I|}\}>0\mbox{ and }\eta_2:=\min_{I\in \mathcal{F}^0}\{\rho_{\textbf{i}_0(I)},\rho_{\textbf{i}_1(I)}\}>0.$$
Write $\eta:=\eta_1\eta_2$ which will be used later. Here $0<\eta\leq 1$.

In the following, we will always assume Assumption A and B.

Under these assumptions, we then have the following another  blow-up
principle.

\textbf{Lemma 4.3.}\emph{ If $J\subset [0,1]$ is any interval of the
form $[0,x]$, then there exists another interval $J'=[0,x']$ of the
same form with $\rho_1<x'\leq 1$ such that $d(J')=d(J)$; similarly,
if $J$ is any interval of the form $[y,1]$, then there exists
another interval $J'=[y',1]$ of the same form with $0\leq
y'<1-\rho_m$ such that $d(J')=d(J)$.}

\emph{Proof.} By Assumption B, if $J\subset S_{1}([0,1])$,
$S_{1}^{-1}J$ is a larger interval of the same density. We iterate
this procedure until we obtain $J'$ not lying in $S_1([0,1])$. The
proof of the second case is similar. $\Box$

For an island $I=[c,d]\subset [0,1]$, we introduce the following notations.
$$\underline{D}_0(I)=\inf_{0<x\leq 1}\{d([c,c+x(d-c)])\} \mbox{ and } \underline{D}_1(I)=\inf_{0<x\leq 1}\{d([d-x(d-c),d])\}.$$
Obviously, if $I_1$ and $I_2$ are two islands with the same overlap
type, i.e., $[I_1]=[I_2]$, then
$$\underline{D}_0(I_1)=\underline{D}_0(I_2) \mbox{ and } \underline{D}_1(I_1)=\underline{D}_1(I_2).$$
Hence the notations $\underline{D}_0(I)$ and $\underline{D}_1(I)$
depend only on the overlap type of the island $I$. Thus we could
define the following constants. For $1\leq i\leq q$, define
$$\underline{D}_0^i:=\underline{D}_0(I) \mbox{ and } \underline{D}_1^i:=\underline{D}_1(I)$$ where $I$ is an arbitrary $\mathcal{T}_i$ type island. These notations are independent of the choice of the island $I$. The following lemma shows Assumption A and Assumption B will ensure that all $\underline{D}_0^i$'s are equal and all $\underline{D}_1^i$'s are also equal.

\textbf{Lemma 4.4.} \emph{\begin{equation}\label{d}
\underline{D}_0^1=\cdots=\underline{D}_0^q \mbox{ and }
\underline{D}_1^1=\cdots=\underline{D}_1^q.
\end{equation}}

\emph{Proof.} Fix $2\leq i\leq q$. We prove $\underline{D}_0^1=\underline{D}_0^i$. Take a $\mathcal{T}_i$ type island $I:=[c,d]$. Now we turn to prove $\underline{D}_0([0,1])=\underline{D}_0(I)$ since $[0,1]$ is of $\mathcal{T}_1$ type.

Using Remark 2, we have $S_{\textbf{v}_0}S_{\textbf{i}_0}([0,1])\cap
K=S_{\textbf{v}_0}S_{\textbf{i}_0} K$ where
$\textbf{v}_0=\textbf{v}_0(I)$ and $\textbf{i}_0=\textbf{i}_0(I)$.
Hence for any $0<x\leq 1$, since
$S_{\textbf{v}_0}S_{\textbf{i}_0}([0,1])\cap K$ is similar to $K$,
we have
$$d([0,x])=d([c,c+x(S_{\textbf{v}_0}S_{\textbf{i}_0}(1)-c)]),$$ which yields that $\underline{D}_0^i\leq \underline{D}_0^1$ by the arbitrariness  of $x$ and the fact that $[c,c+x(S_{\textbf{v}_0}S_{\textbf{i}_0}(1)-c)]\subset I$.

On the other hand, since $I=\bigcup_{\textbf{v}\in
V(I)}S_\textbf{v}([0,1])$, we denote by $c=a_1<a_2<\cdots<a_n$ the
left endpoints of all constitutive intervals of $I$ in increasing
order. (By Assumption A, it is impossible that some distinct
constitutive intervals share a common left endpoint.)  For any
interval $[c,z]\subset I$, choose a largest $a_i$ such that $a_i\leq
z$. If $a_i<z$,
\begin{eqnarray*}
d([c,z])&=&\frac{\lambda([a_1,a_2])+\cdots+\lambda([a_{i-1},a_i])+\lambda([a_i,z])}{((a_2-a_1)+\cdots+(a_i-a_{i-1})+(z-a_i))^\alpha}\\
&\geq&\frac{\lambda([a_1,a_2])+\cdots+\lambda([a_{i-1},a_i])+\lambda([a_i,z])}{(a_2-a_1)^\alpha+\cdots+(a_i-a_{i-1})^\alpha+(z-a_i)^\alpha}\\
&\geq&\min\{d([a_1,a_2]),\cdots,d([a_{i-1},a_i]),d([a_i,z])\}\\
&\geq& \underline{D}_0^1.
\end{eqnarray*}

The last inequality follows from the fact that for each $\textbf{v}\in V(I)$, $S_\textbf{v}K\subset S_\textbf{v}([0,1])\cap K$ and $S_\textbf{v}K$ is similar to $K$. In the case $a_i=z$, we have the same result by a similar discussion.
The arbitrariness of $z$ yields that $\underline{D}_0^i\geq \underline{D}_0^1$.

Hence we have $\underline{D}_0^1=\underline{D}_0^i$. By the arbitrariness of $2\leq i\leq q$, we get the first equality in $(\ref{d})$. The second equality can be proved in a similar way. $\Box$

Thus we could define the common value of
$\underline{D}_0^1,\cdots,\underline{D}_0^q$ as $\underline{D}_0$
and the common value of $\underline{D}_1^1,\cdots,\underline{D}_1^q$
as $\underline{D}_1$. We would like to characterize
$\underline{D}_0$ and $\underline{D}_1$ by the parameters $\beta$'s
and $\gamma$'s.

\textbf{Lemma 4.5.} \emph{Let $k$ be the smallest integer such that
$\beta_1^{(k)}\leq \rho_1$. Then
\begin{equation*}
\underline{D}_0=\min\{d([0,x]): x>0 \mbox{ and }[0,x]\in \mathcal{F}_{k}\}.
\end{equation*}
Similarly, let $k'$ be the smallest integer such that
$\beta_{l_{k'}}^{(k')}\leq \rho_m$. Then
\begin{equation*}
\underline{D}_1=\min\{d([y,1]): y<1 \mbox{ and }[y,1]\in \mathcal{F}_{k'}\}.
\end{equation*}
}
\emph{Proof.} For simplicity, we only prove the first equality. By the definition of $\underline{D}_0$, we note that
$$\underline{D}_0=\inf_{0<x\leq 1}d([0,x]).$$
By Assumption B and the blow-up principle Lemma 4.3, we only need to
consider the interval $[0,x]$ with $\rho_1<x\leq 1$. Since
$d([0,x])$ is a continuous function of $x$, $d([0,x])$ attains its
minimum $\underline{D}_0$ at some $x_0\in [\rho_1,1]$. Furthermore
noting that $d([0,1])=\underline{D}_0$ whenever
$d([0,\rho_1])=\underline{D}_0$, we can assume $x_0>\rho_1$.

It is clear that the point $x_0$ can not fall in a non-empty lake of $\mathcal{F}_{k}$ because then
 $[0,x_0]$ would not have minimal density.
Therefore there exists a $k$-th generation island $I=[c,d]$ such
that $x_0\in [c,d]$. Here $c>0$ because otherwise
$x_0\leq\beta_1^{(k)}\leq \rho_1$ which contradicts $x_0>\rho_1$.
Take $u=x_0-c$.  Assume that $u>0$, then
$$
d([0,x_0])=\frac{\lambda([0,c])+\lambda([c,x_0])}{(c+u)^\alpha}>\frac{\lambda([0,c])+\lambda([c,x_0])}{c^\alpha+u^\alpha}\geq\min\{d([0,c]),
d([c,x_0])\}\geq\underline{D}_0,
$$
which contradicts the minimality of $d([0,x_0])$. Hence $u=0$. Thus
$$\underline{D}_0=\min\{d([0,x]): x>0 \mbox{ and }[0,x]\in \mathcal{F}_{k}\}. \quad \Box$$

We also need to consider the maximal value of the density $d(J)$ whenever $J$ is of the form $[0,x]$ or $[y,1]$. For this purpose, define
$$\overline{D}_0=\sup\{d([0,x]): 0<x\leq 1\} \mbox{ and } \overline{D}_1=\sup\{d([y,1]): 0\leq y<1\}.$$
We would also like to characterize $\overline{D}_0$ and
$\overline{D}_1$ by the parameters $\beta$'s and $\gamma$'s. The
following elementary calculus lemma by Ayer \& Strichartz is useful.

\textbf{Lemma 4.6$^{\cite{Ayer}}$.} \emph{Suppose $0<\alpha<1$,
$p\leq p_0$, $a\geq a_0$, $\kappa>0$ and $y\geq \kappa x^\alpha$.
Then
\begin{equation}\label{a}
0<x\leq \left(\frac{a_0\kappa}{p_0}\right)^{\frac{1}{1-\alpha}}
\end{equation}
implies
\begin{equation}\label{b}
\frac{p-y}{(a-x)^\alpha}<\frac{p}{a^\alpha}.
\end{equation}}

To make this paper self-contained, we give the proof of Lemma 4.6 as
follows.

\emph{Proof of Lemma 4.6.} Consider the function $f(x)=(p-\kappa
x^\alpha)/(a-x)^\alpha$. Noting that $(p-y)/(a-x)^\alpha\leq f(x)$
by the assumption $y\geq \kappa x^\alpha$. And by $p/a^\alpha=f(0)$,
it suffices to show $f'(x)<0$ on the interval $0<x<
(a_0\kappa/p_0)^{1/(1-\alpha)}$. This can be verified by a direct
computation. $\Box$

Using the above lemma, we have the following result concerning $\overline{D}_0$ and $\overline{D}_1$.

\textbf{Lemma 4.7.} \emph{Let $k$ be the smallest integer such that
$ \beta_{\max}^{(k)}\leq (\rho_1 \underline{D}_1)^{1/(1-\alpha)}. $
Then
$$\overline{D}_0=\max\{d([0,x]): x>0 \mbox{ and }[0,x]\in \mathcal{F}_k\}.$$ Similarly, let $k'$ be the smallest integer such that
$ \beta_{\max}^{(k')}\leq (\rho_m \underline{D}_0)^{1/(1-\alpha)}. $
Then
$$\overline{D}_1=\max\{d([y,1]): y<1 \mbox{ and }[y,1]\in \mathcal{F}_{k'}\}.$$
} \emph{Proof.} For simplicity, we only prove the first equality. By
the blow-up principle Lemma 4.3 we can take $x_0> \rho_1$ such that
the interval $[0,x_0]$ has the maximal density (By compactness of
$[\rho_1,1]$, the maximum is attained). The point $x_0$ cannot fall
in a non-empty lake of $\mathcal{F}_k$ because, if so, $[0,x_0]$
would not have maximal density. Let $[0,a]$ be the smallest interval
in $\mathcal{F}_k$ that contains $[0,x_0]$. Then $x_0=a-x$ for some
$x\leq \beta_{\max}^{(k)}$ because $\beta_{\max}^{(k)}$ is the
length of the largest island in $\mathcal{F}^0_k$. Set
$p=\lambda([0,a])$ and $y=\lambda([a-x,a])$, then
$d([0,a-x])=(p-y)/(a-x)^\alpha$ and $d([0,a])=p/a^\alpha$. Thus the
conclusion $(\ref{b})$ of Lemma 4.6 would give $d([0,a-x])<d([0,a])$
unless $x=0$, which implies that $[0,a]$ attains the maximal
density.

To complete the proof we will verify the hypothesis of Lemma 4.6
with $p_0=1$, $a_0=\rho_1$ and $\kappa=\underline{D}_1$. We already
know $a\geq \rho_1$, and $p\leq 1$ is trivial. To verify $y\geq
\kappa x^\alpha$, we observe that $y/x^\alpha=d([a-x,a])$, and by
the definition of  $\underline{D}_1$, we immediately get
$d([a-x,a])\geq \underline{D}_1$. The hypothesis of Lemma 4.6 is
verified, and condition $(\ref{a})$ follows from $x\leq
\beta_{\max}^{k}$ and the hypothesis $ \beta_{\max}^{(k)}\leq
(\rho_1 \underline{D}_1)^{1/(1-\alpha)}. $ $\Box$
\subsection{The maximal density and the Hausdorff measure}

 As
stated before, we denote by $k_0$ the smallest non-negative integer
such that none of the islands in $\mathcal{F}_{k_0+1}^0$ is of new
overlap type. First, in the case that all lakes are non-zero, i.e.,
$\gamma_{\min}^{(k_0+1)}>0$, we have the following result.

\textbf{Theorem 4.8.} \emph{Assume $\gamma_{\min}^{(k_0+1)}>0$, and
let $k\geq k_0+1$ be the smallest integer such that
\begin{equation}\label{e}
2\beta_{\max}^{(k)}\leq (\gamma_{\min}^{(k_0+1)}\min\{{\underline{D}_0,\underline{D}_1}\})^{\frac{1}{1-\alpha}}.
 \end{equation}
 Then the maximal density $d_{\max}$ is attained for an interval in $\mathcal{F}_k$.}

\emph{Proof.} By the blow-up principle Lemma 4.1 we may focus
attention to intervals containing at least one lake of
$\mathcal{F}_{k_0+1}$. So we have the lower bound
$\gamma_{\min}^{(k_0+1)}$ for the length of the interval, which
implies by compactness that the maximal density $d_{\max}$ is
attained. If $[x_1,x_2]$ is an interval of maximal density, we let
$[z_1,z_2]$ be the smallest interval in $\mathcal{F}_k$ containing
$[x_1,x_2]$. Write $a=z_2-z_1$ for the length of the interval,
$x=(z_2-z_1)-(x_2-x_1)$ for the difference of the lengths,
$p=\lambda([z_1,z_2])$ and
$y=\lambda([z_1,x_1])+\lambda([x_2,z_2])$, so that
$d([x_1,x_2])=(p-y)/(a-x)^\alpha$ and $d([z_1,z_2])=p/a^\alpha$. By
a similar way in proving Lemma 4.7, we will complete the proof by
applying Lemma 4.6.

We take $p_0=1$ and $a_0=\gamma_{\min}^{(k_0+1)}$, so that $a\geq a_0$ and $p\leq p_0$. We choose $\kappa=\min\{\underline{D}_0,\underline{D}_1\}$. For the right side interval $[x_2,z_2]$ we have
$$\frac{\lambda([x_2,z_2])}{(z_2-x_2)^\alpha}\geq \underline{D}_1,$$ and similarly for the left side interval $[z_1,x_1]$ we have
$$\frac{\lambda([z_1,x_1])}{(x_1-z_1)^\alpha}\geq \underline{D}_0.$$
Thus we have
$$y\geq \min\{\underline{D}_0,\underline{D}_1\}((x_1-z_1)^\alpha+(z_2-x_2)^\alpha)\geq \kappa x^\alpha.$$ Thus the hypothesis of Lemma 4.6 is verified, and condition $(\ref{a})$ follows from $(\ref{e})$ since $x$ is the sum of two terms, $x_1-z_1$ and $z_2-x_2$, each being at most $\beta_{\max}^{(k)}$. $\Box$

Now we turn to discuss the case that there  exist touching islands.
First we can still obtain a result if we assume a logarithmic
arithmetic relation between $\rho_1$ and $\rho_m$. In the following,
$\eta$ is the positive number defined in Remark 3.

\textbf{Theorem 4.9.} \emph{Suppose there exist positive integers
$n_1$ and $n_m$ such that $\rho_1^{n_1}=\rho_m^{n_m}$. Let $k$ be
the smallest integer such that
\begin{equation*}
2\beta_{\max}^{(k)}\leq(\eta\beta_{\min}^{(k_0+1)}\rho_1^{n_1}\min\{\underline{D}_0,\underline{D}_1\})^{\frac{1}{1-\alpha}}.
\end{equation*}
Then the maximal density $d_{\max}$ is attained for an interval in $\mathcal{F}_k$.}

\textit{Proof.} We claim that it suffices to consider intervals of length at least $\eta \beta_{\min}^{(k_0+1)}\rho_1^{n_1}$. To see this we need a variant of the blow-up principle that shows how to replace smaller intervals with larger intervals of greater density.

Start with any interval $I_0$ not contained in a $(k_0+1)$-th
generation island. If it actually contains a $(k_0+1)$-th generation
island, its length is at least $\beta_{\min}^{(k_0+1)}$, and we are
done. If not, it begins at a point in $I$ and ends at a point in
$I'$ where $I$, $I'$ are two adjacent $(k_0+1)$-th generation
islands. Denote by $L$ the lake separating $I$ and $I'$. Using
Remark 2, there exists a vertex $\textbf{v}\in V(I)$ and an index
$\textbf{i}=(m,\cdots,m)$ such that $S_{\textbf{v}}([0,1])\subset I$
has the same right endpoint as that of $I$ and
$S_{\textbf{v}}S_{\textbf{i}}([0,1])\cap
K=S_{\textbf{v}}S_{\textbf{i}}K$. Similarly,  there also exists a
vertex $\textbf{v}'\in V(I')$ and an index
$\textbf{i}'=(1,\cdots,1)$ such that $S_{\textbf{v}'}([0,1])\subset
I'$ has the same left endpoint as that of $I'$ and
$S_{\textbf{v}'}S_{\textbf{i}'}([0,1])\cap
K=S_{\textbf{v}'}S_{\textbf{i}'}K$. For simplicity, we denote
$S_{\textbf{v}}S_{\textbf{i}}([0,1])$ and
$S_{\textbf{v}'}S_{\textbf{i}'}([0,1])$ by $\widetilde{I}$ and
$\widetilde{I}'$ respectively. (These notations will be used again
in Lemma 4.11. and Theorem 4.12.) Now consider the intervals
$J=S_\textbf{v} S_\textbf{i} S_m^{n_m}([0,1])$ and
$J'=S_{\textbf{v}'}S_{\textbf{i}'}S_1^{n_1}([0,1])$ which lie on the
extreme ends of the lake $L$ separating $\widetilde{I}$ and
$\widetilde{I}'$. These intervals have length
$\rho_\textbf{v}\rho_\textbf{i}\rho_m^{n_m}$ and
$\rho_{\textbf{v}'}\rho_{\textbf{i}'}\rho_{1}^{n_1}$. Moreover,
$J\cap K$ and $J'\cap K$ are similar to $K$. If our interval $I_0$
contains one of them, we are done.

Next suppose our interval begins with a point in $J$ and ends with a point in $J'$, say $I_0=J_0\cup L\cup J'_0$ where $J_0=I_0\cap J$ and $J'_0=I_0\cap J'$ and $L$ is the lake separating $J$ and $J'$. We generate another interval $I_1=J_1\cup L\cup J'_1$ by blowing up $J_0$ to $J_1$ and $J'_0$ to $J'_1$ by a factor $\rho_1^{-n_1}=\rho_m^{-n_m}$. Specifically, we set $J_1=S_\textbf{v}S_\textbf{i}S_m^{-n_m}(S_\textbf{v}S_\textbf{i})^{-1}J_0$ and
$J'_1=S_{\textbf{v}'}S_{\textbf{i}'}S_1^{-n_1}(S_{\textbf{v}'}S_{\textbf{i}'})^{-1}J'_0$. Note that $S_\textbf{v}S_\textbf{i}S_m^{-n_m}(S_\textbf{v}S_\textbf{i})^{-1}$ maps $J$ onto $\widetilde{I}$ and fixes the right endpoint, while $S_{\textbf{v}'}S_{\textbf{i}'}S_1^{-n_1}(S_{\textbf{v}'}S_{\textbf{i}'})^{-1}$ maps $J'$ onto $\widetilde{I}'$ and fixes the left endpoint. So $I_1$ is an interval. We have
$$d(I_0)=\frac{\lambda(J_0)+\lambda(J'_0)}{(|J_0|+|L|+|J'_0|)^\alpha},$$
while
$$
d(I_1)=\frac{\rho_{m}^{-n_m
\alpha}\lambda(J_0)+\rho_{1}^{-n_1\alpha}\lambda(J'_0)}{(\rho_{m}^{-n_m}|J_0|+|L|+\rho_{1}^{-n_1}|J'_0|)^\alpha}=\frac{\lambda(J_0)+\lambda(J'_0)}{(|J_0|+\rho_1^{n_1}|L|+|J'_0|)^\alpha}.
$$
So $d(I_1)\geq d(I_0)$. By iterating this blow-up construction we
eventually arrive at an interval containing either $J$ or $J'$ whose
density is greater than the original interval $I_0$. Hence we only
need to consider intervals containing either $J$ or $J'$ which have
length at least
$\min\{\rho_\textbf{v}\rho_\textbf{i},\rho_{\textbf{v}'}\rho_{\textbf{i}'}\}\rho_1^{n_1}$.
This completes the proof that it suffices to look at intervals of
length at least $\eta\beta_{\min}^{(k_0+1)}\rho_1^{n_1}$.

The rest of the argument is identical to the proof of Theorem 4.8,
except that we take $a_0=\eta\beta_{\min}^{(k_0+1)}\rho_1^{n_1}$.
$\Box$

We consider now the case when the contraction ratios $\rho_1$ and
$\rho_m$ do not satisfy the arithmetic condition. Another elementary
calculus lemma proved by Ayer \& Strichartz in $\cite{Ayer}$ will be
needed.

\textbf{Lemma 4.10$^{{\cite{Ayer}}}$.} \emph{Let $a,a',q,q'>0$, and
$0<\alpha<1$. $F(x)$ is a function
\begin{equation}\label{g}
F(x)=\frac{q+q' x^\alpha}{(a+a' x)^\alpha}
\end{equation}
of positive variables. Then $F$ attains the maximal value of $
\left((\frac{q}{a^\alpha})^{\frac{1}{1-\alpha}}+(\frac{q'}{a'^\alpha})^{\frac{1}{1-\alpha}}\right)^{1-\alpha}
$ at the point $ x_0=(\frac{aq'}{a'q})^{\frac{1}{1-\alpha}}. $
Furthermore, $F(x)$ is strictly increasing on $0\leq x\leq x_0$ and
strictly decreasing on $x>x_0$.}

\emph{Proof.} This can be done by a directly computation of $F'(x)$. $\Box$

\textbf{Lemma  4.11.} \emph{Suppose $I$ and $I'$ are two
$(k_0+1)$-th generation touching islands, and $\rho_1$ and $\rho_m$
are non-arithmetic in the sense that $\rho_1^{n_1}\neq\rho_m^{n_m}$
for any positive integers $n_1$ and $n_m$. $\widetilde{I}$ and
$\widetilde{I}'$ are the same as that defined in the proof of
Theorem 4.9. Then the maximal density of intervals beginning in
$\widetilde{I}$ and ending in $\widetilde{I}'$ is
\begin{equation}\label{j}
(\overline{D}_0^{\frac{1}{1-\alpha}}+\overline{D}_1^{\frac{1}{1-\alpha}})^{1-\alpha}.
\end{equation}}

\emph{Proof.} Any such interval can be written as
$I_0=\widetilde{I}_0\cup \widetilde{I}'_0$ where
$\widetilde{I}_0\subset \widetilde{I}$ ends at the right endpoint of
$\widetilde{I}$, and $\widetilde{I}'_0\subset \widetilde{I}'$ begins
at the left endpoint of $\widetilde{I}'$. Let
$\textbf{v},\textbf{v}',\textbf{i},\textbf{i}'$ be the notations
used in the proof of Theorem 4.9. So
$S_\textbf{v}S_\textbf{i}([0,1])=\widetilde{I}$ and
$S_{\textbf{v}'}S_{\textbf{i}'}([0,1])=\widetilde{I}'$. For any
positive integers $k$ and $k'$ we can form the interval
$$I(k,k')=S_\textbf{v}S_\textbf{i}S_m^{k}(S_\textbf{v}S_\textbf{i})^{-1}\widetilde{I}_0\cup S_{\textbf{v}'}S_{\textbf{i}'}S_1^{k'}(S_{\textbf{v}'}S_{\textbf{i}'})^{-1}\widetilde{I}'_0,$$
which contracts $\widetilde{I}_0$ by a factor of $\rho_m^{k}$ and $\widetilde{I}'_0$ by a factor of $\rho_1^{k'}$, keeping their common endpoint fixed. Then
\begin{equation}\label{k}
d(I(k,k'))=\frac{\rho_m^{k\alpha}\lambda(\widetilde{I}_0)+\rho_1^{k'\alpha}\lambda(\widetilde{I}'_0)}{(\rho_m^{k}|\widetilde{I}_0|+\rho_1^{k'}|\widetilde{I}'_0|)^\alpha}=
\frac{\rho_m^{k\alpha}q+\rho_1^{k'\alpha}q'}{(\rho_m^{k}a+\rho_1^{k'}a')^\alpha},
\end{equation}
where $a=|\widetilde{I}_0|$, $a'=|\widetilde{I}'_0|$ and
$q=\lambda(\widetilde{I}_0)$, $q'=\lambda(\widetilde{I}'_0)$. Notice
that this is exactly of the form $(\ref{g})$ with
$x=\rho_1^{k'}\rho_m^{-k}$, and by the non-arithmetic hypothesis
that $x$ takes on a dense set of values on the positive line. Thus
by Lemma 4.10, $(\ref{k})$ has its maximal value
\begin{equation}\label{l}
(d(\widetilde{I}_0)^{\frac{1}{1-\alpha}}+d(\widetilde{I}'_0)^{\frac{1}{1-\alpha}})^{1-\alpha}.
\end{equation}

Since $(\ref{l})$ is an increasing function of $d(\widetilde{I}_0)$ and $d(\widetilde{I}'_0)$, and by the fact that $\widetilde{I}\cap K$ and $\widetilde{I}'\cap K$ are similar to $K$, it is clear that its maximum is attained when $d(\widetilde{I}_0)$ and $d(\widetilde{I}_0')$ assume their maxima, and these are clearly $\overline{D}_1$ and $\overline{D}_0$. Hence the maximal density of intervals beginning in $\widetilde{I}$ and ending in $\widetilde{I}'$ is $(\overline{D}_0^{1/(1-\alpha)}+\overline{D}_1^{1/(1-\alpha)})^{1-\alpha}$. $\Box$

\textbf{Theorem 4.12.} \emph{Suppose $\gamma_{\min}^{(k_0+1)}=0$,
and $\rho_1$ and $\rho_m$ are non-arithmetic. Let $k$, $k_1$ and
$k_2$ be the smallest integers such that
\begin{equation*}
2\beta_{\max}^{(k)}\leq
(\eta\min\{\beta_{\min}^{(k_0+1)},\gamma_{\min}^{*(k_0+1)}\}\cdot\min\{\underline{D}_0,\underline{D}_1\})^{\frac{1}{1-\alpha}},
\end{equation*}
\begin{equation*}
\beta_{\max}^{(k_1)}\leq
(\rho_1\underline{D}_1)^{\frac{1}{1-\alpha}}\mbox{ and
}\beta_{\max}^{(k_2)}\leq
(\rho_m\underline{D}_0)^{\frac{1}{1-\alpha}}.
\end{equation*}
 Then the maximal density $d_{\max}$ is equal to the maximum of the
finite set of values $d(I)$ as $I$ varies over all intervals in
$\mathcal{F}_{k}$, and
$(d(I_1)^{1/(1-\alpha)}+d(I_2)^{1/(1-\alpha)})^{1-\alpha}$ as $I_1$
varies over all intervals of the form $[0,x]$ in $\mathcal{F}_{k_1}$
and $I_2$ varies over all intervals of the form $[y,1]$ in
$\mathcal{F}_{k_2}$. }

\emph{Proof.}  It follows from $\gamma_{\min}^{(k_0+1)}=0$ that
there exist $(k_0+1)$-th generation touching islands in
$\mathcal{F}^0_{k_0+1}$. If $I$ and $I'$ are two such islands with
$I$ lying on the left side of $I'$, we apply Lemma 4.11, which means
we have to consider the values of $(\ref{j})$. But by Lemma 4.7,
$\overline{D}_0$ is attained for an interval of the form $[0,x]$ in
$\mathcal{F}_{k_1}$, and   $\overline{D}_1$ is attained for an
interval of the form $[y,1]$ in $\mathcal{F}_{k_2}$.

For every two touching $(k_0+1)$-th generation islands $I$ and $I'$,
denote $\widetilde{I}$ and $\widetilde{I}'$ the corresponding
subsets of $I$ and $I'$ respectively as that discussed in Lemma
4.11. We need to consider all intervals  beginning in $I\setminus
\widetilde{I}$ and ending in $I'$, or beginning in $I$ and ending in
$I'\setminus \widetilde{I}'$ for some touching $(k_0+1)$-th
generation islands $I$ and $I'$, and all intervals that contain
either a non-zero lake of $\mathcal{F}_{k_0+1}$ or a $(k_0+1)$-th
generation island. In the first case, intervals have length at least
$\min\{|\widetilde{I}|,|\widetilde{I}'|\}$ which always greater than
$\beta_{\min}^{(k_0+1)}\eta$ by Remark 3. And in the second case,
intervals have lengths bounded below by
$\min\{\beta_{\min}^{(k_0+1)},\gamma_{\min}^{*(k_0+1)}\}$. Hence the
lengths of all the above two kinds of intervals are greater than
$\eta\min\{\beta_{\min}^{(k_0+1)},\gamma_{\min}^{*(k_0+1)}\}$. Then
by a slightly variant of the proof of Theorem 4.8, the maximal
density over all such intervals is attained by an interval in
$\mathcal{F}_k. \Box$

\subsection{The minimal centered density and the packing measure}

First, we give a lemma concerning the relation between $d_{\min}$ and $\underline{D}_0,\underline{D}_1$.

\textbf{Lemma 4.13.} \emph{$d_{\min}\leq
2^{-\alpha}\min\{\underline{D}_0, \underline{D}_1\}$.}

\emph{Proof.} It suffices to show that there exist $J_0$, $J_1$
centered in $K$, contained in $[0,1]$, such that
$d(J_0)=2^{-\alpha}\underline{D}_0$ and
$d(J_1)=2^{-\alpha}\underline{D}_1$. For simplicity, we only prove
the first equality. By Lemma 4.5, let $k$ be the smallest integer
such that $\beta_1^{(k)}\leq \rho_1$, then there exists $x_0$ with
$[0,x_0]\subset \mathcal{F}_k$ such that
$d([0,x_0])=\underline{D}_0$. Since $\alpha<1$, there must exist at
least one non-empty lake in $\mathcal{F}_1$, i.e., there exists
$1\leq i\leq {l_1}-1$ with $\gamma_{i}^{(1)}>0$, where
$\gamma_i^{(1)}$ is the length of the lake separating the $i$-th and
$(i+1)$-th first generation islands $I$ and $I'$. From the
discussion in Remark 2, there is a unique vertex $\textbf{v}\in
V(I')$ and an index $\textbf{i}=(1,\cdots,1)$ such that the
sub-interval $S_\textbf{v}S_\textbf{i}([0,1])\subset I'$ has the
same left endpoint as that of $I'$ and furthermore,
$S_\textbf{v}S_\textbf{i}([0,1])\cap K=S_\textbf{v}S_\textbf{i} K$,
i.e., $S_\textbf{v}S_\textbf{i}([0,1])\cap K$ is similar to $K$ with
a contraction ratio $\rho_\textbf{v}\rho_\textbf{i}$. Choose a
non-negative integer $k'$ large enough so that
$\rho_\textbf{v}\rho_\textbf{i}\rho_1^{k'}x_0<\gamma_i^{(1)}$. Hence
the interval
$[S_\textbf{v}S_\textbf{i}S_{1}^{k'}(0)-\rho_\textbf{v}\rho_\textbf{i}\rho_1^{k'}x_0,S_\textbf{v}S_\textbf{i}S_{1}^{k'}(0)]$
is contained in the lake separating $I$ and $I'$. Thus
$\lambda([S_\textbf{v}S_\textbf{i}S_{1}^{k'}(0)-\rho_\textbf{v}\rho_\textbf{i}\rho_1^{k'}x_0,S_\textbf{v}S_\textbf{i}S_{1}^{k'}(0)])=0$.
Define
$$J_0:=[S_\textbf{v}S_\textbf{i}S_{1}^{k'}(0)-\rho_\textbf{v}\rho_\textbf{i}\rho_1^{k'}x_0,S_\textbf{v}S_\textbf{i}S_{1}^{k'}(0)+\rho_\textbf{v}\rho_\textbf{i}\rho_1^{k'}x_0].$$
Since the interval $S_\textbf{v}S_\textbf{i}S_1^{k'}([0,1])\cap K$
is also similar to $K$, we have
\begin{eqnarray*}
d(J_0)&=&\frac{\lambda([S_\textbf{v}S_\textbf{i}S_{1}^{k'}(0)-\rho_\textbf{v}\rho_\textbf{i}\rho_1^{k'}x_0,S_\textbf{v}S_\textbf{i}S_{1}^{k'}(0)+\rho_\textbf{v}\rho_\textbf{i}\rho_1^{k'}x_0])}{2^\alpha (\rho_\textbf{v}\rho_\textbf{i}\rho_1^{k'}x_0)^\alpha}\\
&=&\frac{\lambda(S_\textbf{v}S_\textbf{i}S_1^{k'}([0,x_0]))}{2^\alpha
|S_\textbf{v}S_\textbf{i}S_1^{k'}([0,x_0])|^\alpha}=
2^{-\alpha}d(S_\textbf{v}S_\textbf{i}S_1^{k'}([0,x_0]))\\&=&2^{-\alpha}d([0,x_0])=2^{-\alpha}\underline{D}_0,
\end{eqnarray*}
which concludes the proof. $\Box$

\textbf{Lemma 4.14.} \emph{Let $J\subset [0,1]$ be an interval
centered in $K$ and not contained in any $(k_0+1)$-th generation
island. If
$|J|\leq\min\{\beta_{\min}^{(k_0+1)},\gamma_{\min}^{*(k_0+1)}\},$
then $d(J)\geq 2^{-\alpha}\min\{\underline{D}_0,\underline{D}_1\}.$}

\emph{Proof.} For each interval $J\subset [0,1]$ centered in $K$ and not contained in any $(k_0+1)$-th generation island with $|J|\leq \min\{\beta_{\min}^{(k_0+1)},\gamma_{\min}^{*(k_0+1)}\}$. It is clear that there are only three possible cases for $J$.

Case 1: There exist two touching $(k_0+1)$-th generation islands $I_1$ and $I_2$ with $I_1$ lying on the left side of $I_2$, and $J=J_1\cup J_2$, where $J_1\subset I_1$ and $J_2\subset I_2$.

In this case, we have
$$
d(J)=\frac{\lambda(J_1)+\lambda(J_2)}{(|J_1|+|J_2|)^\alpha}
>\frac{\lambda(J_1)+\lambda(J_2)}{|J_1|^\alpha+|J_2|^\alpha}
\geq \min\{d(J_1),d(J_2)\}
\geq \min\{\underline{D}_0,\underline{D}_1\}
$$
by the definition of the constants $\underline{D}_0,
\underline{D}_1$. It follows that
$d(J)>\min\{\underline{D}_0,\underline{D}_1\}$.

Case 2: There exist two separate $(k_0+1)$-th generation islands
$I_1$ and $I_2$ with $I_1$ lying on the left side of $I_2$, and
$J=J_1\cup J_2$, where $J_1\subset I_1$, $J_2\subset L$, and $L$ is
the lake separating $I_1$ and $I_2$.

In this case, we have $|J_1|\geq |J_2|$ since $J$ is centered in $K$. Thus
$$d(J)=\frac{\lambda(J_1)}{(|J_1|+|J_2|)^\alpha}\geq\frac{\lambda(J_1)}{2^\alpha|J_1|^\alpha}=2^{-\alpha}d(J_1)\geq 2^{-\alpha}\underline{D}_1$$
by the definition of $\underline{D}_1$.

Case 3: There exist two separate $(k_0+1)$-th generation islands
$I_1$ and $I_2$ with $I_1$ lying on the left side of $I_2$, and
$J=J_1\cup J_2$, where $J_1\subset L$, $J_2\subset I_2$, and $L$ is
the lake separating $I_1$ and $I_2$.

In this case, we have $d(J)\geq 2^{-\alpha}\underline{D}_0$ by a discussion similar to Case 2.

Combining the above discussion completes the proof. $\Box$

For convenience, denote all $(k_0+1)$-th generation islands in increasing order by $I_1,\cdots,I_{l_{k_0+1}}$. For each $1\leq i\leq l_{k_0+1}$, write $I_{i}=[a_i,b_i]$.

\textbf{Lemma 4.15.} \emph{If
$d_{\min}<2^{-\alpha}\min\{\underline{D}_0,\underline{D}_1\}$, then
$$d_{\min}=\min_{1\leq i_1<i_2<l_{k_0+1}}\frac{\sum_{i=i_1+1}^{i_2}\lambda([a_i,b_i])}{(a_{i_2+1}-b_{i_1}-2\mbox{dist}(\frac{a_{i_2+1}+b_{i_1}}{2},K))^\alpha}.$$}

\emph{Proof.} Since
$d_{\min}<2^{-\alpha}\min\{\underline{D}_0,\underline{D}_1\}$, by
the blow-up principle Lemma 4.1 and Lemma 4.14, we only need to
consider
 intervals  with lengths greater than $\min\{\beta_{\min}^{(k_0+1)},\gamma_{\min}^{*(k_0+1)}\}$ which are centered in $K$ and not contained in any $(k_0+1)$-th generation island.

By the compactness of $K$, there exists a such interval $J_0=[a_0,b_0]$ such that $d_{\min}=d(J_0).$  First we prove the following statements.

(1) Either $a_0\in\{b_i: 1\leq i\leq l_{k_0+1}-1\}$ or $a_0$ is contained in one lake;

(2) Either $b_0\in\{a_i: 2\leq i\leq l_{k_0+1}\}$ or $b_0$ is contained in one lake.

For simplicity we only prove $(1)$. The statement $(2)$ will follow by a similar argument. Assume that $(1)$ is not true. Then there exists a $1\leq i\leq l_{k_0+1}-1$ such that $a_0\in [a_i,b_i)$. In the following we will lead to a contradiction. We first claim $(a_0+b_0)/2>b_i$. Otherwise $(a_0+b_0)/2\in [a_i,b_i]$, then
$$d([a_0,b_0])=\frac{\lambda([a_0,b_0])}{(b_0-a_0)^\alpha}\geq\frac{\lambda([a_0,b_i])}{2^\alpha(b_i-a_0)^\alpha}=2^{-\alpha}d([a_0,b_i])\geq 2^{-\alpha}\underline{D}_1,$$
which contradicts the assumption $d_{\min}<2^{-\alpha}\min\{\underline{D}_0,\underline{D}_1\}$. Hence it follows that
\begin{eqnarray*}
d([a_0,b_0])&=&\frac{\lambda([a_0,b_i])+\lambda([b_i,a_0+b_0-b_i])+\lambda([a_0+b_0-b_i,b_0])}{(2(b_i-a_0)+(a_0+b_0-2b_i))^\alpha}\\
&>&\frac{\lambda([a_0,b_i])+\lambda([b_i,a_0+b_0-b_i])+\lambda([a_0+b_0-b_i,b_0])}{2^\alpha(b_i-a_0)^\alpha+(a_0+b_0-2b_i)^\alpha}\\
&\geq&\frac{\lambda([a_0,b_i])+\lambda([b_i,a_0+b_0-b_i])}{2^\alpha(b_i-a_0)^\alpha+(a_0+b_0-2b_i)^\alpha}\\
&\geq&\min\{2^{-\alpha}\frac{\lambda([a_0,b_i])}{(b_i-a_0)^\alpha},\frac{\lambda([b_i,a_0+b_0-b_i])}{(a_0+b_0-2b_i)^\alpha}\}\\
&=&\min\{2^{-\alpha}d([a_0,b_i]),d([b_i,a_0+b_0-b_i])\}\\
&\geq&\min\{2^{-\alpha}\underline{D}_1,d([b_i,a_0+b_0-b_i])\}.
\end{eqnarray*}
Since $[b_i,a_0+b_0-b_i]\subset [0,1]$ is an interval centered in
$K$, the above inequality contradicts the fact that $d([a_0,b_0])$
attains the minimal value $d_{\min}<2^{-\alpha} \underline{D}_1$.
Thus the statement $(1)$ is true.

By the statement $(1)$ and $(2)$, we have $ a_0\in
[b_{i_1},a_{i_1+1}) \mbox{ and } b_0\in (b_{i_2},a_{i_2+1}] $ for
some $1\leq i_1<i_2<l_{k_0+1}$. Since
$\lambda([a_0,b_0])=\sum_{i=i_1+1}^{i_2}\lambda([a_i,b_i])$ and
$d([a_0,b_0])$ attains the minimal centered density $d_{\min}$, it
follows that $a_0$, $b_0$ are taken such that $(b_0-a_0)$ is the
largest value under the condition $a_0\in[b_{i_1},a_{i_1+1})$,
$b_{0}\in(b_{i_2},a_{i_2+1}]$ and $(a_0+b_0)/2\in K$. Thus we have
$b_0-a_0=a_{i_2+1}-b_{i_1}-2
\mbox{dist}(\frac{a_{i_2+1}+b_{i_1}}{2},K)$ and
$$d([a_0,b_0])=\frac{\sum_{i=i_1+1}^{i_2}\lambda([a_i,b_i])}{(a_{i_2+1}-b_{i_1}-2 \mbox{dist}(\frac{a_{i_2+1}+b_{i_1}}{2},K))^\alpha}.$$

Therefore we complete the proof of Lemma 4.15. $\Box$

\textbf{Theorem 4.16.} $d_{\min}=\min\{2^{-\alpha}\underline{D}_0,
2^{-\alpha}\underline{D}_1,D\}.$ Here
$$\underline{D}_0=\min\{d([0,x]): x>0 \mbox{ and }[0,x]\in \mathcal{F}_{k}\},$$ where $k$ is the smallest number such that $\beta_1^{(k)}\leq \rho_1$;
$$
\underline{D}_1=\min\{d([y,1]): y<1 \mbox{ and }[y,1]\in \mathcal{F}_{k'}\},
$$ where $k'$ is the smallest number such that $\beta_{l_{k'}}^{(k')}\leq \rho_m$;
$$D=\min_{1\leq i_1<i_2<l_{k_0+1}}\frac{\sum_{i=i_1+1}^{i_2}\lambda([a_i,b_i])}{(a_{i_2+1}-b_{i_1}-2\mbox{dist}(\frac{{a_{i_2+1}+b_{i_1}}}{2},K))^\alpha},$$
where $[a_i,b_i]$ is the $i$-th $(k_0+1)$-th generation island for each $1\leq i\leq l_{k_0+1}$.

\emph{Proof.} It follows immediately from Lemma 4.5, Lemma 4.13 and
Lemma 4.15. $\Box$

\subsection{Examples on computing measures}

In this subsection we will show how to compute the Hausdorff and
packing measures of $K$, using the above theory. Let's look at the
examples in Section 2 again.

\textbf{Example 4.17.} If $\{S_j=\rho_j x+b_j\}_{j=1}^m$ satisfies
the OSC as showed in Example 2.4. We assume $0<\rho_j<1$ for each
$j=1,\cdots,m$. Without loss of generality, we assume that the
images $S_j([0,1])$ are in increasing order, with $S_1(0)=0$ and
$S_m(1)=1$. Then $\mathcal{H}^\alpha(K)$ and $\mathcal{P}^\alpha(K)$
can be calculated as the minimal or maximal value of a finite set of
elementary functions of the parameters $\rho$'s and $b$'s. This is
true since $K$ naturally satisfies all the assumptions in this
section. The results have already been proved in $\cite{Ayer}$ and
$\cite{Feng}$ for $\mathcal{H}^\alpha(K)$ and
$\mathcal{P}^\alpha(K)$ respectively.

\textbf{Example 4.18.} Let $K$ be the invariant set of the IFS
$\{S_j\}_{j=1}^3$ on $\mathbb{R}$ defined in Example 2.5. Let
$\alpha$ denotes the dimension of $K$. If $\rho+2r-\rho r=1$, then
$\alpha=1$, $K=[0,1]$, thus $\mathcal{H}^1(K)=\mathcal{P}^1(K)=1$.
Hence we only need to consider the non-trivial case $\rho+2r-\rho
r<1$. In this case,  $\alpha<1$, $k_0$=1. Using Lemma 4.5, we get
$$\underline{D}_0=\min\{d([0,x]): x>0 \mbox{ and } [0,x]\in \mathcal{F}_2\}\mbox{ and }\underline{D}_1=\min\{d([y,1]): y<1 \mbox{ and } [y,1]\in \mathcal{F}_1\}.$$ We adopt the same notations of Example 2.5, with $\mathcal{T}_1=[I_1]$ and $\mathcal{T}_2=[I_2]$
 where $I_1=[0,1]$ is the root island and $I_2=S_{1}([0,1])\cup S_2([0,1])$. Using the definition of $\lambda$, we have $\lambda(I_1)=1$ and $\lambda(I_2)=1-r^\alpha$.
  Then after a detailed computation, we get the exact values $$\underline{D}_0=\min\{1,\frac{1-r^\alpha}{(1-r)^\alpha},\frac{(\rho^\alpha+r^\alpha)(1-r^\alpha)}{(\rho+r)^\alpha(1-r)^\alpha},\frac{1-r^{2\alpha}}{(1-r^2)^\alpha}\}=\frac{1-r^\alpha}{(1-r)^\alpha},$$
   and $$\underline{D}_1=\min\{1,\frac{r^\alpha}{(1-\rho-r+\rho r)^\alpha}\}=\frac{r^\alpha}{(1-\rho-r+\rho r)^\alpha}.$$

\emph{Hausdorff measure.} A detailed calculation shows
$\gamma_{\min}^{(k_0+1)}=\gamma_{\min}^{(2)}=(1-2r-\rho+\rho
r)\cdot\min\{\rho,r\},$ and $\beta_{\max}^{(k)}=(\rho+r-\rho
r)\cdot(\max\{\rho,r\})^{k-1}$ for each $k\geq 1$.
 Let $k\geq 2$ be the smallest integer
such that $ 2\beta_{\max}^{(k)}\leq
(\gamma_{\min}^{(2)}\min\{{\underline{D}_0,\underline{D}_1}\})^{\frac{1}{1-\alpha}},
$ i.e.,
\begin{equation*}
2(\rho+r-\rho r)\cdot(\max\{\rho,r\})^{k-1}\leq ((1-2r-\rho+\rho
r)\cdot\min\{\rho,r\}\cdot\min\{{\underline{D}_0,\underline{D}_1}\})^{\frac{1}{1-\alpha}}.
 \end{equation*}
 Then by Theorem 4.8, the maximal density $d_{\max}$ is attained for an interval in $\mathcal{F}_k$. Furthermore, by Corollary 1.3, $\mathcal{H}^\alpha(K)=d_{\max}^{-1}$.

\emph{Packing measure.} Since $k_0=1$, the constant $D$ in Theorem
4.16 is
$$D=\min_{1\leq i_1<i_2<5}\frac{\sum_{i=i_1+1}^{i_2}\lambda([a_i,b_i])}{(a_{i_2+1}-b_{i_1}-2\mbox{dist}(\frac{a_{i_2+1}+b_{i_1}}{2},K))^\alpha}.$$
where $[a_i,b_i]$ is the $i$-th $2$-th generation island for each
$1\leq i\leq 5$. A detailed calculation yields the exact value of
$D$,
\begin{eqnarray*}
D=\min\{&&\frac{r^\alpha-r^{2\alpha}}{(\rho+r-2\rho r-r^2-\rho^2+\rho^2 r-2\mbox{dist}(\frac{\rho+r-r^2+\rho^2-\rho^2 r}{2},K))^\alpha},\\
&&\frac{r^\alpha}{(1-r-\rho^2-\rho r+\rho^2 r-2\mbox{dist}(\frac{1-r+\rho^2+\rho r-\rho^2 r}{2},K))^\alpha},\\
&&\frac{2r^\alpha-r^{2\alpha}}{(1-r^2-\rho^2-\rho r+\rho^2 r-2\mbox{dist}(\frac{1-r^2+\rho^2+\rho r-\rho^2 r}{2},K))^\alpha},\\
&&\frac{r^{2\alpha}}{(1-r-\rho-r^2+\rho r^2-2\mbox{dist}(\frac{1-r+\rho+r^2-\rho r^2}{2},K))^\alpha},\\
&&\frac{r^\alpha}{(1-\rho-2r^2+\rho r^2-2\mbox{dist}(\frac{1+\rho-\rho r^2}{2},K))^\alpha},\\
&&\frac{r^\alpha-r^{2\alpha}}{(1-r-\rho-r^2+\rho
r-2\mbox{dist}(\frac{1+r+\rho-r^2-\rho r}{2},K))^\alpha}\}.
\end{eqnarray*}
Then by Theorem 4.16, the minimal centered density
 $d_{\min}=\min\{2^{-\alpha}\underline{D}_0,2^{-\alpha}\underline{D}_1,D\}.$ Hence by Corollary 1.4, $\mathcal{P}^\alpha(K)=d_{\min}^{-1}$. $\Box$

\textbf{Remark.} Consider the special case where $\rho=r=1/16$ in
the above example. By Example 2.5, the dimension
$\alpha=\log_{16}{2}/({3-\sqrt{5}})\approx 0.3471.$ And we calculate
that $\underline{D}_0=(16^\alpha-1)/{15^\alpha}\approx 0.6320$,
$\underline{D}_{1}=16^\alpha/225^\alpha\approx0.3995$, and
$\gamma_{\min}^{(2)}=209/{4096}\approx0.0510$. Moreover, for each
$k\geq 1$, $\beta_{\max}^{(k)}=({31}/{256})\cdot({1}/{16^{k-1}})$.
Hence the smallest $k$ should satisfy
$$2\cdot\frac{31}{256}\cdot \frac{1}{16^{k-1}}\leq(0.0510\cdot 0.3995)^{\frac{1}{1-\alpha}},$$
which yields that $k=3$. After a complicated computation using
computer, we eventually get
$d_{\max}=({256^{\alpha}-16^{\alpha}})/{31^\alpha}\approx1.2861$,
and
$\mathcal{H}^{\alpha}(K)={31^\alpha}/({256^\alpha-16^\alpha})\approx
0.7775.$

For the packing measure, a detailed calculation using computer
shows that $D={8^\alpha}/{225^\alpha}\approx 0.3140.$ Hence
$d_{\min}=2^{-\alpha}\underline{D}_1=D={8^\alpha}/{225^\alpha}\approx
0.3140$ and
$\mathcal{P}^\alpha(K)={225^\alpha}/{8^\alpha}\approx3.1843.$

 \textbf{Example 4.19.} Let $K$ be the invariant set of the IFS $\{S_j\}_{j=1}^3$ on $\mathbb{R}$ defined in Example 2.6.
Then the dimension $\alpha$ of $K$ is the logarithmic ratio of the
largest root of the polynomial equation $x^3-6x^2+5x-1=0$ to $9$,
$\alpha\approx0.7369$.

 We adopt the same notations of
Example 2.6, with $\mathcal{T}_1=[I_1]$, $\mathcal{T}_2=[I_2]$ and
$\mathcal{T}_3=[I_3]$ where $I_1=[0,1]$, $I_2=S_{11}([0,1])\cup
S_{12}([0,1])$ and $I_3=S_{13}([0,1])\cup S_2([0,1])$. Using the
definition of $\lambda$, we have $\lambda(I_1)=1$,
$\lambda(I_2)=1/3^\alpha-1/9^\alpha$ and
$\lambda{(I_3)}=1-2/3^\alpha+1/9^\alpha$. Using Lemma 4.5, we get
$$\underline{D}_0=\min\{d([0,x]): x>0 \mbox{ and } [0,x]\in \mathcal{F}_1\}\mbox{ and }\underline{D}_1=\min\{d([y,1]): y<1 \mbox{ and } [y,1]\in \mathcal{F}_1\}.$$ A detailed calculation shows that
$\underline{D}_0=\min\{1,\frac{3^\alpha-1}{2^\alpha},\frac{9^\alpha-1}{8^\alpha}\}=\frac{3^\alpha-1}{2^\alpha}\approx0.7482,$
and $\underline{D}_1=\min\{1,\frac{9^\alpha}{16^\alpha},
\frac{81^\alpha-27^\alpha+9^\alpha}{70^\alpha}\}=\frac{9^\alpha}{16^\alpha}\approx0.6544.$

\emph{Hausdorff measure.} Since $k_0=1$, a detailed calculation
shows
$\gamma_{\min}^{(k_0+1)}=\gamma_{\min}^{(2)}=\frac{7}{729}\approx0.0096.$
By Theorem 4.8, we need to find a smallest integer $k$  such that
\begin{eqnarray*}
2\beta_{\max}^{(k)}\leq
(\gamma_{\min}^{(2)}\min\{{\underline{D}_0,\underline{D}_1}\})^{\frac{1}{1-\alpha}}=(\frac{7}{729}\min\{\frac{3^\alpha-1}{2^\alpha},\frac{9^\alpha}{16^\alpha}\})^{\frac{1}{1-\alpha}}.
 \end{eqnarray*}
 Noticing that $\beta_{\max}^{(k)}=5/(27\cdot 9^{k-1})$, the smallest $k=10$.

Then by Theorem 4.8, the maximal density $d_{\max}$ is attained for
an interval in $\mathcal{F}_{10}$. After a complicated computation
by using computer, we eventually get
$d_{\max}=(27^{\alpha}-9^{\alpha})/{11^\alpha}\approx1.0756$, and
$\mathcal{H}^{\alpha}(K)={11^\alpha}/({27^\alpha-9^\alpha})\approx
0.9297.$

 \emph{Packing measure.} Since $k_0=1$, the constant $D$ in Theorem 4.16 is
$$D=\min_{1\leq i_1<i_2<20}\frac{\sum_{i=i_1+1}^{i_2}\lambda([a_i,b_i])}{(a_{i_2+1}-b_{i_1}-2\mbox{dist}(\frac{a_{i_2+1}+b_{i_1}}{2},K))^\alpha}.$$
where $[a_i,b_i]$ is the $i$-th $2$-th generation island for each
$1\leq i\leq 20$. A detailed calculation by using computer yields
the exact value of $D=9^\alpha/32^\alpha\approx0.3927$. Hence
$d_{\min}=2^{-\alpha}\underline{D}_1=D={9^\alpha}/{32^\alpha}\approx
0.3927$ and
$\mathcal{P}^\alpha(K)={32^\alpha}/{9^\alpha}\approx2.5467.$ $\Box$

Although these algorithms applies in theory to any case considered
under Assumption A and B, in practice it is useable in very few
cases. Even in the following simple example.

\textbf{Example 4.20.} Let $K$ be the invariant set of the IFS
$\{S_j\}_{j=1}^4$ on $\mathbb{R}$ defined in Example 2.7. The
dimension $\alpha$ of $K$ equals
${\log_4(5+\sqrt{5})}-\frac{1}{2}\approx0.9276$.

 We adopt the same and notations of Example 2.7, with
$\mathcal{T}_1=[I_1]$, $\mathcal{T}_2=[I_2]$ and
$\mathcal{T}_3=[I_3]$ where $I_1=[0,1]$, $I_2=S_{2}([0,1])\cup
S_{3}([0,1])$ and $I_3=S_{22}([0,1])\cup S_{23}([0,1])\cup
S_{31}([0,1])$. Using the definition of $\lambda$, we have
$\lambda(I_1)=1$, $\lambda(I_2)=1-2/4^\alpha$ and
$\lambda{(I_3)}=1-3/4^\alpha$. Using Lemma 4.5, we get
$$\underline{D}_0=\min\{d([0,x]): x>0 \mbox{ and } [0,x]\in \mathcal{F}_1\}\mbox{ and }\underline{D}_1=\min\{d([y,1]): y<1 \mbox{ and } [y,1]\in \mathcal{F}_1\}.$$ A detailed calculation shows that
$\underline{D}_0=\min\{1,\frac{4^\alpha-1}{3^\alpha}\}=\frac{4^\alpha-1}{3^\alpha}\approx0.9449,$
and $\underline{D}_1=\min\{1,\frac{2^\alpha}{3^\alpha},
\frac{4^\alpha-1}{3^\alpha}\}=\frac{2^\alpha}{3^\alpha}\approx0.6865.$

\emph{Hausdorff measure.} We will use Theorem 4.9 since there exist
touching islands and $\rho_1=\rho_m$. Observe that $k_0=2$,
$n_1=n_m=1$ and $\eta=1/8$. A detailed calculation shows $
\beta_{\min}^{(k_0+1)}=\beta_{\min}^{(3)}=\frac{1}{64}.$ By Theorem
4.9, we need to find a smallest integer $k$  such that
$$
2\beta_{\max}^{(k)}\leq(\eta\beta_{\min}^{(k_0+1)}\rho_1^{n_1}\min\{\underline{D}_0,\underline{D}_1\})^{\frac{1}{1-\alpha}}=(\frac{1}{8}\cdot\frac{1}{64}\cdot\frac{1}{4}\min\{\frac{4^\alpha-1}{3^\alpha},\frac{2^\alpha}{3^\alpha}\})^{\frac{1}{1-\alpha}}.
$$
 Noticing that $\beta_{\max}^{(k)}=1/(8\cdot 4^{k-2})=1/2^{2k-1}$, the smallest $k=81$.

 Then by Theorem 4.9, the maximal density $d_{\max}$ is attained for an interval in $\mathcal{F}_{81}$. Furthermore, by Corollary 1.3, $\mathcal{H}^\alpha(K)=d_{\max}^{-1}$.
However, the time involved in searching all sets in
$\mathcal{F}_{k}$ rapidly becomes impractical. Hence our algorithm
for computing $\mathcal{H}^\alpha(K)$ exceeds the computing power.

 \emph{Packing measure.} Since $k_0=2$, the constant $D$ in Theorem 4.16 is
$$D=\min_{1\leq i_1<i_2<35}\frac{\sum_{i=i_1+1}^{i_2}\lambda([a_i,b_i])}{(a_{i_2+1}-b_{i_1}-2\mbox{dist}(\frac{a_{i_2+1}+b_{i_1}}{2},K))^\alpha}.$$
where $[a_i,b_i]$ is the $i$-th $3$-th generation island for each
$1\leq i\leq 35$. A detailed calculation by using computer yields
the exact value of $D={1}/{3^\alpha}\approx0.3609$. Hence
$d_{\min}=2^{-\alpha}\underline{D}_1=D={1}/{3^\alpha}\approx0.3609$
and $\mathcal{P}^\alpha(K)=3^\alpha\approx2.7706.$ $\Box$

\section{Further Discussions}

\textbf{Are Assumption A and B necessary?}

If we permit the IFS not to satisfy Assumption A or  B, things
become more complicated. It seems hard to get a general formulae  of
$\mathcal{H}^\alpha(K)$ and $\mathcal{P}^\alpha(K)$. However, in
some special cases, we can still use the similar method to get the
results. The following are two concrete examples.

\textbf{Example 5.1.} Consider the IFS $\{S_j\}_{j=1}^3$ as follows.
$$S_1(x)=\frac{1}{3}x,\quad S_2(x)=\frac{1}{9}x+\frac{2}{9},\quad S_3(x)=\frac{1}{3}x+\frac{2}{3}.$$
If we choose $\mathcal{M}_k=\Lambda_{k}$ for each $k\geq 0$, then
$\{S_j\}_{j=1}^3$ satisfies the GTFC with respect to the invariant
set $(0,1)$, with Assumption A not satisfied. See Figure 4.
\setlength{\unitlength}{0.5cm}
\begin{figure}[htbp]
\begin{center}
\includegraphics[width=14cm]{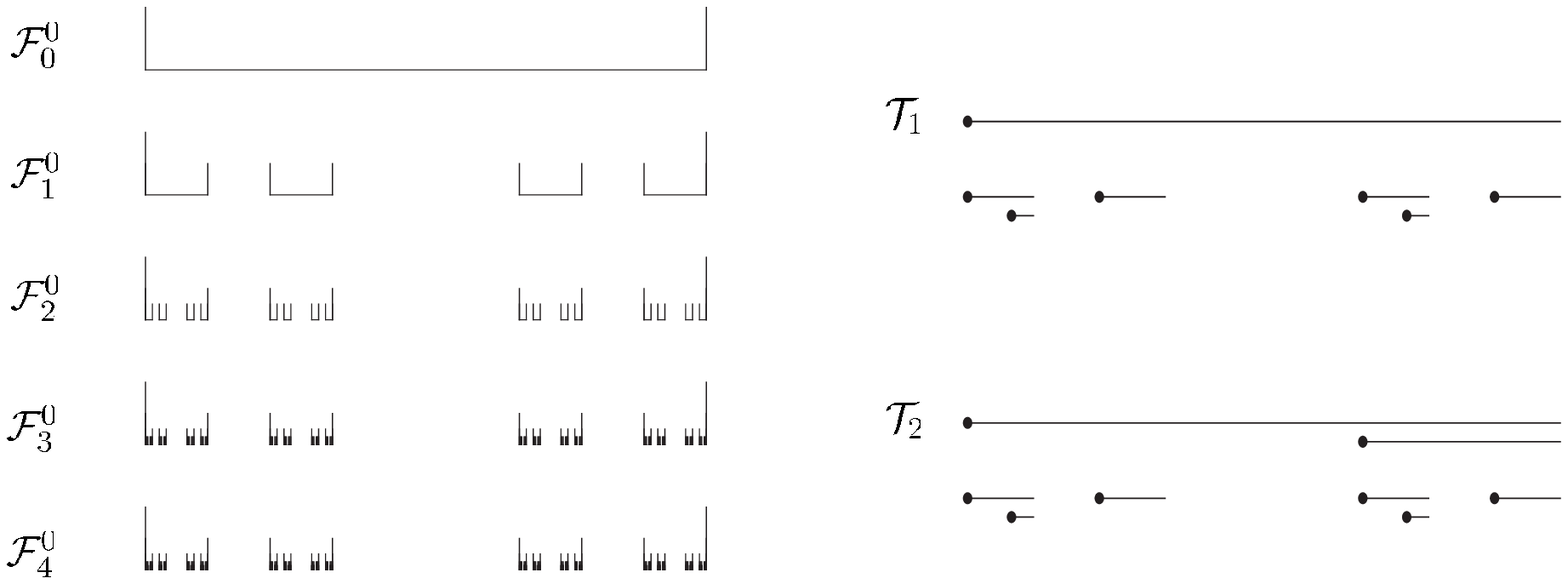}\\
\end{center}\emph{ Figure 4.  The
first five levels of islands and  the distinct overlap types in
Example 5.1.}
\end{figure}
In fact, the invariant set $K$ of this IFS is the classical Cantor
ternary set. It is well-known that the  dimension of $K$ is equal to
$\log_32$ and $\mathcal{H}^{\log_32}=1$. Replacing the IFS by
$\{S_1,S_3\}$ equivalently, from Theorem 4.16, one gets that
$\mathcal{P}^{\log_32}=4^{\log_32}$. (See this also in
$\cite{Feng}$.)

\textbf{Example 5.2.} Consider the IFS $\{S_j\}_{j=1}^3$ defined as
follows.
$$S_1(x)=\rho x,\quad S_2(x)=\rho x+\rho\frac{1-\rho}{1+\rho},\quad S_3(x)=\rho
x+1-\rho,$$ where $0<\rho<1/3$. Choose $\mathcal{M}_k=\Sigma_{k}$
for each $k\geq 0$, then $\{S_j\}_{j=1}^3$ satisfies the GTFC with
respect to the invariant set $(0,1)$. We omit the proof. Observing
that $\rho(1-\rho)/(1+\rho)\in S_1([0,1])\cap K$ and
$\rho(1-\rho)/(1+\rho)\notin S_1K$, one gets that Assumption B is
not satisfied. By Theorem 2.3, it is easy to verify that the
dimension of the $K$ is $\alpha=\log_{\rho}2/(3+\sqrt{5})$. See
Figure 5.

It is easy to verify that $[I_1]$ and $[I_2]$ denoted respectively
by $\mathcal{T}_1$ and $\mathcal{T}_2$ are the all distinct overlap
types, where $I_1=[0,1]$ and $I_2=S_1([0,1])\cup S_2([0,1])$. In
this case Lemma 4.4 may not hold since Assumption B is not
satisfied. Hence the original notations $\underline{D}_0$ and
$\underline{D}_1$ may not be suitable. However, we redefine them by
$\underline{D}_0:=\min\{\underline{D}_0^1,\underline{D}_0^2\}$, and
$\underline{D}_1:=\min\{\underline{D}_1^1,\underline{D}_1^2\}$. A
similar discussion as the proof of Lemma 4.4 shows that
$\underline{D}_0^1\leq\underline{D}_0^2$ and
$\underline{D}_1^1\leq\underline{D}_1^2$. Moreover, observing the
distribution of the offsprings of $I_2$, we can also get
$\underline{D}_1^2\leq\underline{D}_1^1$. Hence
$\underline{D}_0=\underline{D}_0^1$ and
$\underline{D}_1=\underline{D}_1^1=\underline{D}_1^2$. By a similar
argument of Lemma 4.5 (with suitable modifications), we get
$$\underline{D}_0=\min\{d([0,x]): x>0 \mbox{ and } [0,x]\in
\mathcal{F}_2\}\mbox{ and }\underline{D}_1=\min\{d([y,1]): y<1
\mbox{ and } [y,1]\in \mathcal{F}_1\}.$$ Using the definition of
$\lambda$, we have $\lambda(I_1)=1$ and
$\lambda(I_2)=1-\rho^\alpha$. Hence we can get that $d([0,x])$
attains the minimal value
$\underline{D}_0=\min\{\frac{2(1-\rho^\alpha)(1+\rho)^\alpha}{(1-\rho)^\alpha(2+\rho)^\alpha},\frac{1-\rho^\alpha}{(1-\rho)^\alpha}\}=\frac{1-\rho^\alpha}{(1-\rho)^\alpha}$
at the point $x_0=1-\rho$, and $d([y,1])$ attains the minimal value
$\underline{D}_1=\rho^\alpha\frac{(1+\rho)^\alpha}{(1-\rho)^\alpha}$
at the point $y_0=2\rho/(1+\rho)$.

\emph{Hausdorff measure.} $k_0=1$,
$\gamma_{\min}^{(2)}=\rho(1-3\rho)/(1+\rho)$,
$\beta_{\max}^{(k)}=2\rho^k/(1+\rho)$ for each $k\geq 1$. By a
suitable modification of Theorem 4.8, we get $d_{\max}$ is attained
for an interval in $\mathcal{F}_k$ where $k\geq 2$ is the smallest
integer such that
$2\beta_{\max}^{(k)}\leq(\gamma_{\min}^{(2)}\min\{\underline{D}_0,\underline{D}_1\})^{\frac{1}{1-\alpha}},$
namely,
$$
{4\rho^k}\leq
(\frac{(1-3\rho)\rho^{1+\alpha}}{(1-\rho)^\alpha})^{\frac{1}{1-\alpha}}.
$$
 Furthermore, by Corollary 1.3,
$\mathcal{H}^\alpha(K)=d_{\max}^{-1}$.

\emph{Packing measure.} We need a similar result of Lemma 4.13,
i.e., $d_{\min}\leq
2^{-\alpha}\min\{\underline{D}_0,\underline{D}_1\}$. However, at
first glance, we can not prove it in general. The reason is the
following. Recall the proof of Lemma 4.13, we should find two
intervals $J_0$ and $J_1$ centered in $K$ with
$d(J_0)=2^{-\alpha}\underline{D}_0$ and
$d(J_1)=2^{-\alpha}\underline{D}_1$. In fact, we can define $J_1$ as
that used in the proof of Lemma 4.13 since
$\underline{D}_1^1=\underline{D}_1^2$  and $S_3([0,1])\cap K=S_3K$.
(Half of Assumption B holds.) But for $J_0$, the original process is
invalid since $\underline{D}_0^1$ and $\underline{D}_0^2$ may not be
equal. Fortunately, a detailed verifying shows that the inequality
$\underline{D}_1\leq \underline{D}_0$ always holds, which ensures
that it is useless to find $J_0$, i.e., the existence of $J_1$ is
enough. Hence Lemma 4.13 remains true. Once Lemma 4.13 is proved,
Lemma 4.14, Lemma 4.15 and eventually Theorem 4.16 follow
automatically, in which $\min\{\underline{D}_0,\underline{D}_1\}$
are all replaced by $\underline{D}_1$.

Hence $\mathcal{P}^\alpha(K)=d_{\min}^{-1}$ where
 $d_{\min}=\min\{\frac{\rho^\alpha(1+\rho)^\alpha}{2^\alpha(1-\rho)^\alpha},D\},$
 and
\begin{eqnarray*}
D=\min\{&&\frac{\rho^\alpha-\rho^{2\alpha}}{(\frac{2\rho-3\rho^2-\rho^3}{1+\rho}-2\mbox{dist}(\frac{2\rho+\rho^2-\rho^3}{2(1+\rho)},K))^\alpha},
\frac{\rho^\alpha}{(\frac{1-3\rho^2}{1+\rho}-2\mbox{dist}(\frac{1+\rho^2}{2(1+\rho)},K))^\alpha},\\
&&\frac{2\rho^\alpha-\rho^{2\alpha}}{(\frac{1+\rho-3\rho^2-\rho^3}{1+\rho}-2\mbox{dist}(\frac{1+\rho+\rho^2-\rho^3}{2(1+\rho)},K))^\alpha},
\frac{\rho^{2\alpha}}{(1-2\rho-2\mbox{dist}(\frac{1}{2},K))^\alpha},\\
&&\frac{\rho^\alpha}{(1-\rho-\rho^2-2\mbox{dist}(\frac{1+\rho-\rho^2}{2},K))^\alpha},
\frac{\rho^\alpha-\rho^{2\alpha}}{(\frac{1-\rho-\rho^2-\rho^3}{1+\rho}-2\mbox{dist}(\frac{1+3\rho-\rho^2-\rho^3}{2(1+\rho)},K))^\alpha}\}.
\Box
\end{eqnarray*}

\setlength{\unitlength}{0.5cm}
\begin{figure}[htbp]
\begin{center}
\includegraphics[width=14cm]{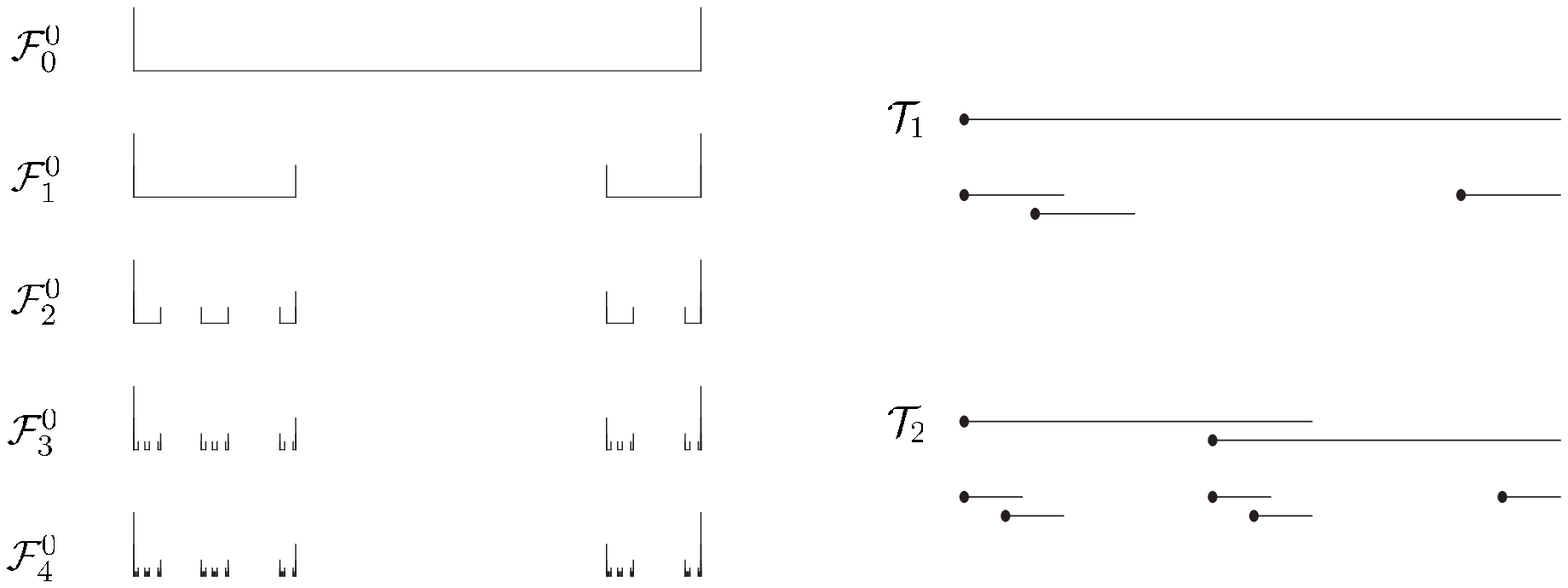}\\
\end{center}\emph{ Figure 5.  The
first five levels of islands and  the distinct overlap types in
Example 5.2.}
\end{figure}

\textbf{Remark.} Consider the special case where $\rho=1/6$ in the
above example. The dimension
$\alpha=\log_{6}({(3+\sqrt{5})}/{2})\approx 0.5371.$ And we
calculate that $\underline{D}_0=(6^\alpha-1)/{5^\alpha}\approx
0.6816$, $\underline{D}_{1}=7^\alpha/30^\alpha\approx0.4576$, and
$\gamma_{\min}^{(2)}=1/{14}\approx0.0714$. Moreover, for each $k\geq
1$, $\beta_{\max}^{(k)}=({2}/{7})\cdot({1}/{6^{k-1}})$. Hence the
smallest $k$ should satisfy
$$2\cdot\frac{2}{7}\cdot \frac{1}{6^{k-1}}\leq(0.0714\cdot 0.4576)^{\frac{1}{1-\alpha}},$$
which yields that $k=5$. After a complicated computation by using
computer, we eventually get
$d_{\max}=2({6^{\alpha}-1})/{6^\alpha}\approx1.2361$, and
$\mathcal{H}^{\alpha}(K)={6^\alpha}/(2(6^\alpha-1))\approx 0.8090.$

For the packing measure, a detailed calculation by using computer
shows that $D={7^\alpha}/{60^\alpha}\approx 0.3154.$ Hence
$d_{\min}=2^{-\alpha}\underline{D}_1=D={7^\alpha}/{60^\alpha}\approx
0.3154$ and
$\mathcal{P}^\alpha(K)={60^\alpha}/{7^\alpha}\approx3.1709.$

\textbf{Can we allow negative $\rho$'s?}

Consider the IFSs containing orientation reversing similarities,
which simply means that we allow some of the $\rho_j$ to be
negative. In some special cases, it can be proved by using a similar
method that the results of Section 4 still remain true.

We again normalize by assuming that $S_j([0,1])$ are in increasing
order with $S_{1}([0,1])$ containing $0$ and $S_m([0,1])$ containing
$1$. For example, $\rho_1\rho_m>0$, $ \bigcup_{2\leq j\leq m}
S_j([0,1])\cap S_1([0,1])=\emptyset$, $\bigcup_{1\leq j\leq
m-1}S_j([0,1])\cap S_m([0,1])=\emptyset$, and
$\gamma_{\min}^{(k_0+1)}>0$. Obviously, in this case Assumption B is
naturally satisfied. We can still obtain measure results following
the idea of Section 4.

To illustrate this, we assume that $\rho_1$ and $\rho_m$ are both
positive since otherwise we can replace the IFS with its iterated
square, i.e., all compositions $S_iS_j$. Similar to Example 5.2,
Lemma 4.10 may not hold. Hence the original notations
$\underline{D}_0$ and $\underline{D}_1$ will be meaningless. In
spite of this, the following constant defined by
$\underline{D}=\min\{\underline{D}_0^1,\underline{D}_1^1\}$ will
replace the important role of the original notation
$\min\{\underline{D}_0,\underline{D}_1\}$ in Section 4. It is not
hard to verify that $\underline{D}\leq \underline{D}_0^i$ and
$\underline{D}\leq \underline{D}_1^i$ for all $1\leq i\leq q$, where
$q$ is the cardinality of all the distinct overlap types. Moreover,
we could also characterize the constant $\underline{D}$ using a
suitable modification of Lemma 4.5, i.e.,
$$\underline{D}=\min\{d(J): J\in \mathcal{F}_1 \mbox{ and } J \mbox{ is of the form } [0,x]\mbox{ or }[y,1]\}.$$
It is not hard to verify that Theorem 4.8 for Hausdorff measure, and
Lemma 4.11, Lemma 4.14, Lemma 4.15 and Theorem 4.16 for packing
measure still remain true in which
$\min\{\underline{D}_0,\underline{D}_1\}$ is replaced by
$\underline{D}$.

\textbf{How about self-similar sets in higher dimensional Euclidean
spaces?}

As mentioned earlier, with suitable modifications if necessary, the
results in Section 3 may be generalized to self-similar sets in
higher dimensional Euclidean spaces. How about measure results in
Section 4? Our answer is: almost all the obvious generalizations of
our results are false.

It is clear that the blow-up principles continue to hold. So we can
also focus attention to sets not contained entirely in some
$(k_0+1)$-generation island, where $k_0$ is also the smallest
non-negative integer such that none of the islands in
$\mathcal{F}_{k_0+1}$ is of a new overlap type. For the maximal
density $d_{\max}$, since any set can be replaced by its convex hull
without decreasing the density, it is reasonable to limit any
searching algorithm to convex sets. However, it does not follow that
the maximal density among $\mathcal{F}_k$ sets is achieved by a
convex set, since the convex hull of a set in $\mathcal{F}_k$ may
not belong to $\mathcal{F}_k$. To illustrate this, in $\cite{Ayer}$
Ayer \& Strichartz consider a concrete example, i.e, the usual
Sierpinski gasket in plane. We omit it here. For the minimal
centered density $d_{\min}$, one obvious obstacle is that almost all
lemmas concerning $d_{\min}$ require $\alpha<1$. Moreover, even if
we were to limit attention to self-similar sets of dimension
$\alpha<1$, it is unlikely that the same results would hold. In
fact, we would have to confront how to describe $\underline{D}_0$
and $\underline{D}_1$ and how to overcome the difficulty of the
calculation of densities of higher dimensional sets. (It seems
uncontrollable.)

\renewcommand{\baselinestretch}{0.8}

\end{document}